\documentclass[10pt]{amsart}
\usepackage{amssymb}
\usepackage{amscd}
\usepackage{verbatim}
\usepackage{epsfig}

\begin{document}

\newcommand\nin{\notin}
\newcommand\identity{\operatorname{id}}
\newcommand\id{\operatorname{id}}
\newcommand\Id{\operatorname{Id}}
\newcommand\pos{\Real^+}
\newcommand\Rnp{\Real\setminus\{0\}}
\newcommand\nzero{\setminus\{0\}}
\newcommand\Cx{\mathbb{C}}
\newcommand\Cxp{\Cx^+}
\newcommand\Cxm{\Cx^-}
\newcommand\Nat{\mathbb{N}}
\newcommand\halfNat{{\frac{1}{2}}\mathbb{N}}
\newcommand\intgr{\mathbb{Z}}
\newcommand\im{\operatorname{Im}}
\newcommand\re{\operatorname{Re}}
\newcommand\sign{\operatorname{sign}}
\newcommand\codim{\operatorname{codim}}
\newcommand\End{\operatorname{End}}
\newcommand\Ker{\operatorname{Ker}}
\newcommand\Hom{\operatorname{Hom}}
\newcommand\ideal{{\mathcal I}}
\newcommand\Span{\operatorname{span}}
\newcommand\Range{\operatorname{Ran}}
\newcommand\graph{\operatorname{graph}}
\newcommand\slim{\operatornamewithlimits{s-lim}}
\newcommand\diag{\operatorname{diag}}
\newcommand\Rn{\Real^n}
\newcommand\Rm{\Real^m}
\newcommand\RN{\Real^N}
\newcommand\RtN{\Real^{2N}}
\newcommand\RM{\Real^M}
\newcommand\sphere{\mathbb{S}}
\newcommand\Sn{\sphere^{n-1}}
\newcommand\Sm{\sphere^{m-1}}
\newcommand\Snp{\sphere^n_+}
\newcommand\Smp{\sphere^m_+}
\newcommand\SN{\sphere^{N-1}}
\newcommand\SNp{\sphere^N_+}
\newcommand\circlep{\sphere^1_+}
\newcommand\Phom{P_{h}}
\newcommand\Shom{S_{h}}
\newcommand\distance{\operatorname{dist}}
\newcommand\cl{\operatorname{cl}}
\newcommand\interior{\operatorname{int}}
\newcommand\Fa{\operatorname{Fa}}
\newcommand\ff{\operatorname{ff}}
\newcommand\mf{\operatorname{mf}}
\newcommand\cf{\operatorname{cf}}
\newcommand\scf{\operatorname{sf}}
\newcommand\lf{\operatorname{lb}}
\newcommand\rf{\operatorname{rb}}
\newcommand\indfam{{\mathcal K}}
\newcommand\calX{{\mathcal X}}
\newcommand\calK{{\mathcal K}}
\newcommand\calF{{\mathcal F}}
\newcommand\calO{{\mathcal O}}
\newcommand\calC{{\mathcal C}}
\newcommand\calCL{{\mathcal C}_{\text L}}
\newcommand\calCR{{\mathcal C}_{\text R}}
\newcommand\Cinf{\CI}
\newcommand\dist{{\mathcal C}^{-\infty}}
\newcommand\dCinf{\dot\Cinf}
\newcommand\ddist{\dot\dist}
\newcommand\Cj{{\mathcal C}^j}
\newcommand\Linf{L^{\infty}}
\newcommand\bcon{{\mathcal A}}
\newcommand\bconc{{\mathcal A}_{\text{phg}}}
\newcommand\Sch{{\mathcal S}}
\newcommand\temp{\Sch^{\prime}}
\newcommand\Diff{\operatorname{Diff}}
\newcommand\Diffb{\operatorname{Diff}_{\text{b}}}
\newcommand\Diffc{\operatorname{Diff}_{\text{c}}}
\newcommand\Diffsc{\operatorname{Diff}_{\text{sc}}}
\newcommand\DiffI{\operatorname{Diff}_{\text{I}}}
\newcommand\DiffIq{\operatorname{Diff}_{\text{I},q}}
\newcommand\supp{\operatorname{supp}}
\newcommand\ssupp{\operatorname{sing\ supp}}
\newcommand\csupp{\operatorname{cone\ supp}}
\newcommand\esupp{\operatorname{ess\ supp}}
\newcommand\Fr{{\mathcal F}}
\newcommand\Frinv{\Fr^{-1}}
\newcommand\bop{{\mathcal B}}
\newcommand\spec{\operatorname{spec}}
\newcommand\pspec{\spec_{pp}}
\newcommand\cspec{\spec_{c}}
\newcommand\FIO{{\mathcal I}}
\newcommand\SP{\operatorname{SP}}
\newcommand\Symc{S_c}
\newcommand\Symca{S_c^{\alpha}}
\newcommand\Symczero{S_c^{0,...,0}}
\newcommand\sci{{}^{\text{sc}}}
\newcommand\sct{\sci T^*}
\newcommand\scct{\sci\bar{T}^*}
\newcommand\Csc{C_{\text{sc}}}
\newcommand\SNpscd{(\SNp)^2_{\text{sc}}}
\newcommand\scdiag{\Delta_{\text{sc}}}
\newcommand\projscl{\pi^L_{\text{sc}}}
\newcommand\projscr{\pi^R_{\text{sc}}}
\newcommand\scHL{\sci H^{2,0}_{|\zeta|^2-\lambda^2}}
\newcommand\scHrg{\sci H^{2,0}_{\sqrt{g}}}
\newcommand\Hsc{H_{\text{sc}}}
\newcommand\WF{\operatorname{WF}}
\newcommand\WFp{\operatorname{WF^{\prime}}}
\newcommand\WFsc{\operatorname{WF}_{\text{sc}}}
\newcommand\WFscp{\operatorname{WF_{sc}^{\prime}}}
\newcommand\WFC{\operatorname{WF}_C}
\newcommand\WFCi{\operatorname{WF}_{C_i}}
\newcommand\elliptic{\operatorname{ell}}
\newcommand\Psop{\operatorname{\Psi}}
\newcommand\Psiscrs{\operatorname{\Psi_{sc}^{-2,\infty}}}
\newcommand\Psiscr{\operatorname{\Psi_{sc}^{-2,0}}}
\newcommand\Psiscrm{\operatorname{\Psi_{sc}^{0,2}}}
\newcommand\PsiscHam{\operatorname{\Psi_{sc}^{2,0}}}
\newcommand\Psisci{\operatorname{\Psi_{sc}^{*,*}}}
\newcommand\Psiscid{\operatorname{\Psi_{sc}^{0,0}}}
\newcommand\Psiscis{\operatorname{\Psi_{sc}^{0,\infty}}}
\newcommand\Psiscsi{\operatorname{\Psi_{sc}^{-\infty,0}}}
\newcommand\Psiscs{\operatorname{\Psi_{sc}^{-\infty,\infty}}}
\newcommand\Psiscalg{\operatorname{\Psi_{sc}^{\infty,-\infty}}}
\newcommand\nullHam{{\mathcal N}}
\newcommand\charD{\Sigma_{\Delta-\lambda^2}}
\newcommand\charLap{\Sigma_{\Delta-\lambda}}
\newcommand\Snl{\Sn_{\lambda}}
\newcommand\SNl{\SN_{\lambda}}
\newcommand\gammat{\tilde\gamma}
\newcommand\gammasc{\gamma}
\newcommand\taut{\tilde\tau}
\newcommand\PlV{P_V(\lambda)}
\newcommand\PlVc{P_V^{\flat}(\lambda)}
\newcommand\Pl{P_0(\lambda)}
\newcommand\SVl{S_V(\lambda)}
\newcommand\Sjr{S_j^{\reduced}}
\newcommand\Rkp{{\mathcal R}^k_+}
\newcommand\Rkm{{\mathcal R}^k_-}
\newcommand\Rkpm{{\mathcal R}^k_{\pm}}
\newcommand\Phys{{\mathcal P}}
\newcommand\Pc{\overline{\mathcal P}}
\newcommand\pip{\pi^{\perp}}
\newcommand\xit{\tilde\xi}
\newcommand\zetat{\tilde\zeta}
\newcommand\etat{\tilde\eta}
\newcommand\sigmat{\tilde\sigma}
\newcommand\sigmahat{\hat\sigma}
\newcommand\thetat{\tilde\theta}
\newcommand\psit{\tilde\psi}
\newcommand\phit{\tilde\phi}
\newcommand\chit{\tilde\chi}
\newcommand\rhot{\tilde\rho}
\newcommand\xib{\bar\xi}
\newcommand\zetab{\bar\zeta}
\newcommand\thetab{\bar\theta}
\newcommand\etab{\bar\eta}
\newcommand\iotal{\iota_{\lambda}}
\newcommand\rhoat{\rhot_{\alpha_1}}
\newcommand\Lambdat{\tilde\Lambda}
\newcommand\poles{\Lambda'}
\newcommand\rpoles{\Lambda_p}
\newcommand\thresholds{\Lambda}
\newcommand\Vt{\tilde V}
\newcommand\half{{\frac{1}{2}}}
\newcommand\sigmah{\sigma^{1/2}}
\newcommand\bX{\partial X}
\newcommand\Deltabt{\tilde\Delta_0}
\newcommand\strip{\Omega_T}
\newcommand\Kf{K^{\flat}}
\newcommand\Gs{G^{\sharp}}
\newcommand\Gt{\tilde G}
\newcommand\Gb{\bar G}
\newcommand\Osc{\sci\Omega}
\newcommand\Osch{\sci\Omega^{\half}}
\newcommand\Oscmh{\sci\Omega^{-\half}}
\newcommand\Isc{I_{sc}}
\newcommand\Qzl{Q^0_{-\lambda}}
\newcommand\Lie{{\mathcal L}}
\newcommand\bl{{\text b}}
\newcommand\scl{{\text{sc}}}
\newcommand\sccl{{\text{scc}}}
\newcommand\Scl{{\text{3sc}}}
\newcommand\ScLl{{\text{Sc,L}}}
\newcommand\ScRl{{\text{Sc,R}}}
\newcommand\Sccl{{\text{Scc}}}
\newcommand\sfl{{\operatorname{s\Phi}}}
\newcommand\sfi{{}^\sfl}
\newcommand\sus{{\text{sus}}}
\newcommand\Osfh{\sfi\Omega^{\half}}
\newcommand\Isf{I_\sfl}
\newcommand\XXb{X^2_\bl}
\newcommand\XXsc{X^2_\scl}
\newcommand\XXSc{X^2_\Scl}
\newcommand\XXScL{X^2_\ScLl}
\newcommand\XXScR{X^2_\ScRl}
\newcommand\MMsc{M^2_\scl}
\newcommand\Deltab{\Delta_\bl}
\newcommand\Deltasc{\Delta_\scl}
\newcommand\DeltaSc{\Delta_\Scl}
\newcommand\DeltaScL{\Delta_\ScLl}
\newcommand\DeltaScR{\Delta_\ScRl}
\newcommand\prs{\sigma}
\newcommand\Nsc{N_\scl}
\newcommand\Nscp{N_{\scl,p}}
\newcommand\Nff{N_{\ff}}
\newcommand\Nffz{N_{\ff,0}}
\newcommand\Nffzp{N_{\ff,0,p}}
\newcommand\Nffl{N_{\ff,l}}
\newcommand\Nffml{N_{\ff,-l}}
\newcommand\Nmf{N_{\mf}}
\newcommand\Nmfz{N_{\mf,0}}
\newcommand\Nmfl{N_{\mf,l}}
\newcommand\Nmfml{N_{\mf,-l}}
\newcommand\ffb{\operatorname{bf}}
\newcommand\Ffb{\operatorname{bf'}}
\newcommand\ffsc{\operatorname{sf}}
\newcommand\ffSc{\operatorname{sf_C}}
\newcommand\Ffsc{\operatorname{sf'}}
\newcommand\rff{\rho_{\ff}}
\newcommand\rmf{\rho_{\mf}}
\newcommand\rffb{\rho_{\ffb}}
\newcommand\rffsc{\rho_{\ffsc}}
\newcommand\rFfsc{\rho_{\Ffsc}}
\newcommand\rffSc{\rho_{\ffSc}}
\newcommand\rinf{\rho_{\infty}}
\newcommand\CL{C_L}
\newcommand\CR{C_R}
\newcommand\betab{\beta_\bl}
\newcommand\betasc{\beta_\scl}
\newcommand\betaSc{\beta_\Scl}
\newcommand\BetaSc{\bar\beta_\Scl}
\newcommand\betaScL{\beta_\ScLl}
\newcommand\betaScR{\beta_\ScRl}
\newcommand\sfT{{}^\sfl T^*}
\newcommand\sfN{{}^\sfl N^*}
\newcommand\ScT{{}^\Scl T^*}
\newcommand\SccT{{}^\Scl \bar T^*}
\newcommand\ScS{{}^\Scl S^*}
\newcommand\scS{{}^\scl S^*}
\newcommand\Tb{{}^\bl T}
\newcommand\Tsc{{}^\scl T}
\newcommand\Tsf{{}^\sfl T}
\newcommand\TSc{{}^\Scl T}
\newcommand\CSc{C_\Scl}
\newcommand\Lambdasc{{}^\scl\Lambda}
\newcommand\XXXb{X^3_\bl}
\newcommand\XXXsc{X^3_\scl}
\newcommand\XXXSc{X^3_\Scl}
\newcommand\XXXScO{X^3_{\Scl,O}}
\newcommand\XXXScF{X^3_{\Scl,F}}
\newcommand\XXXScS{X^3_{\Scl,S}}
\newcommand\XXXScC{X^3_{\Scl,C}}
\newcommand\KDsc{\operatorname{KD^{\half}_\scl}}
\newcommand\SDsc{\operatorname{SD^{\half}_\scl}}
\newcommand\SDsf{\operatorname{SD^{\half}_\sfl}}
\newcommand\KDSc{\operatorname{KD^{\half}_\Scl}}
\newcommand\KDScEF{\operatorname{KD^{E,F}_\Scl}}
\newcommand\Oh{\operatorname{\Omega^{\half}}}
\newcommand\WFSc{\WF_\Scl}
\newcommand\WFScmf{\WF_{\Scl,\mf}}
\newcommand\WFScff{\WF_{\Scl,\ff}}
\newcommand\WFScs{\WF_{\Scl,\prs}}
\newcommand\WFScp{\WF'_\Scl}
\newcommand\WFScmfp{\WF'_{\Scl,\mf}}
\newcommand\WFScffp{\WF'_{\Scl,\ff}}
\newcommand\WFScsp{\WF'_{\Scl,\prs}}
\newcommand\Diffscc{\Diff_\sccl}
\newcommand\DiffSc{\Diff_\Scl}
\newcommand\Diffsf{\Diff_\sfl}
\newcommand\DiffScc{\Diff_\Sccl}
\newcommand\DiffscI{\Diff_{\scl,\text{I}}}
\newcommand\VscI{\Vf_{\scl,\text{I}}}
\newcommand\DiffsV{\operatorname{Diff}_{\sus(V)}}
\newcommand\DiffsVsc{\operatorname{Diff}_{\sus(V),\scl}}
\newcommand\DiffsVCsc{\operatorname{Diff}_{\sus(V)-C,\scl}}   
\newcommand\Psisc{\Psop_\scl}
\newcommand\Psiscc{\Psop_\sccl}
\newcommand\PsiSc{\Psop_\Scl}
\newcommand\PsiScc{\Psop_\Sccl}
\newcommand\PsiSccml{\Psop^{m,l}_\Sccl}
\newcommand\PsiScxx{\Psop^{*,*}_\Scl}
\newcommand\PsiScml{\Psop^{m,l}_\Scl}
\newcommand\PsiScmz{\Psop^{m,0}_\Scl}
\newcommand\PsiScmmz{\Psop^{-m,0}_\Scl}
\newcommand\PsiSckz{\Psop^{k,0}_\Scl}
\newcommand\PsiScmmml{\Psop^{-m,-l}_\Scl}
\newcommand\Psiscmkk{\Psop^{-k,k}_\scl}
\newcommand\Psiscmmmkk{\Psop^{-m-k,k}_\scl}
\newcommand\Psiscmoo{\Psop^{-1,1}_\scl}
\newcommand\Psiscmz{\Psop^{m,0}_\scl}
\newcommand\Psiscmmz{\Psop^{-m,0}_\scl}
\newcommand\PsiSckmkl{\Psop^{km,kl}_\Scl}
\newcommand\PsiScmplp{\Psop^{m',l'}_\Scl}
\newcommand\PsiScmmpllp{\Psop^{m+m',l+l'}_\Scl}
\newcommand\Psiscml{\Psop^{m,l}_\scl}
\newcommand\PsiScid{\Psop^{0,0}_\Scl}
\newcommand\PsiSczo{\Psop^{0,1}_\Scl}
\newcommand\PsiScmii{\Psop^{-\infty,\infty}_\Scl}
\newcommand\PsiScmiz{\Psop^{-\infty,0}_\Scl}
\newcommand\PsiScmoo{\Psop^{-1,1}_\Scl}
\newcommand\PsisCid{\Psop^{0,0}_{\scl-C}}
\newcommand\PsisC{\Psop_{\scl-C}}
\newcommand\Psiinf{\Psop_{\infty}}
\newcommand\Psiinfid{\Psop_{\infty}^0}
\newcommand\PsiFinf{\Psop_{\infty-\Fr}}
\newcommand\PsisVscml{\Psop^{m,l}_{\sus(V),\scl}}
\newcommand\PsisVsc{\Psop_{\sus(V),\scl}}
\newcommand\PsisVpsc{\Psop_{\sus(V_p),\scl}}
\newcommand\PsisVCSc{\Psop_{\sus(V)-C,\scl}}
\newcommand\SFinf{S_{\infty-\Fr}}
\newcommand\YsVC{Y^2_{\sus(V)-C,\scl}}
\newcommand\ffYsc{\ffsc_{\sus(V)}}
\newcommand\SXC{S(X;C)}
\newcommand\Ios{I_{\text{os}}}
\newcommand\pbL{\pi^2_{\bl,{\text L}}}
\newcommand\pbR{\pi^2_{\bl,{\text R}}}
\newcommand\pscL{\pi^2_{\scl,{\text L}}}
\newcommand\pscR{\pi^2_{\scl,{\text R}}}
\newcommand\PbO{\pi^3_{\bl,{\text O}}}
\newcommand\PscO{\pi^3_{\scl,{\text O}}}
\newcommand\PScO{\pi^3_{\Scl,{\text O}}}
\newcommand\PScF{\pi^3_{\Scl,{\text F}}}
\newcommand\PScC{\pi^3_{\Scl,{\text C}}}
\newcommand\PScS{\pi^3_{\Scl,{\text S}}}
\newcommand\pScL{\pi^2_{\Scl,{\text L}}}
\newcommand\pScR{\pi^2_{\Scl,{\text R}}}
\newcommand\CLF{\CL^F}
\newcommand\CLO{\CL^O}
\newcommand\CLS{\CL^S}
\newcommand\CLC{\CL^C}
\newcommand\DeltaYb{\Delta_{\bl,Y}}
\newcommand\DeltaYsc{\Delta_{\sus-\scl}}
\newcommand\Vf{{\mathcal V}}
\newcommand\Vb{{\mathcal V}_{\bl}}
\newcommand\Vsc{{\mathcal V}_{\scl}}
\newcommand\VSc{{\mathcal V}_{\Scl}}
\newcommand\Vsf{{\mathcal V}_{\sfl}}
\newcommand\VfI{\Vf_{\text{I}}}
\newcommand\VfIq{\Vf_{\text{I},q}}
\newcommand\scH{{}^\scl H}
\newcommand\scHg{\scH_g}
\newcommand\xh{\hat x}
\newcommand\xb{\bar x}
\newcommand\Yh{\hat Y}
\newcommand\Yb{\bar Y}
\newcommand\hb{\bar h}
\newcommand\xih{\hat\xi}
\newcommand\etah{\hat\eta}
\newcommand\muh{\hat\mu}
\newcommand\mut{\tilde\mu}
\newcommand\mub{\bar\mu}
\newcommand\mubh{\widehat{\bar\mu}}
\newcommand\yb{\bar y}
\newcommand\phib{\bar \phi}
\newcommand\ub{\bar u}
\newcommand\Qb{\bar Q}
\newcommand\Wbp{{\bar W}^\perp}
\newcommand\Wp{W^\perp}
\newcommand\Kt{\tilde K}
\newcommand\Wt{\tilde W}
\newcommand\Ut{\tilde U}
\newcommand\xt{\tilde x}
\newcommand\yt{\tilde y}
\newcommand\ft{\tilde f}
\newcommand\fs{f^{\sharp}}
\newcommand\at{\tilde a}
\newcommand\htil{\tilde h}
\newcommand\gt{\tilde g}
\newcommand\Ht{\tilde H}
\newcommand\Mt{\tilde M}
\newcommand\St{\tilde S}
\newcommand\Pt{\tilde P}
\newcommand\Rt{\tilde R}
\newcommand\qt{\tilde q}
\newcommand\Qt{\tilde Q}
\newcommand\Xb{\bar X}
\newcommand\lambdat{\tilde\lambda}
\newcommand\epst{\tilde\epsilon}
\newcommand\At{\tilde A}
\newcommand\Ah{\hat A}
\newcommand\Bh{\hat B}
\newcommand\Gh{\hat G}
\newcommand\Hh{\hat H}
\newcommand\Qh{\hat Q}
\newcommand\Ph{\hat P}
\newcommand\Nh{\hat N}
\newcommand\Sh{\hat S}
\newcommand\Gcal{{\mathcal G}}
\newcommand\GcalC{{\mathcal G}_C}
\newcommand\Jcal{{\mathcal J}}
\newcommand\JcalC{{\mathcal J}_C}
\newcommand\Miff{\ \text{iff}\ }
\newcommand\Mif{\ \text{if}\ }
\newcommand\Mand{\ \text{and}\ }
\newcommand\Mor{\ \text{or}\ }
\newcommand\Mst{\ \text{s.t.}\ }

\newcommand\pa{\partial}
\newcommand\sub{\operatorname{sub}}
\newcommand\fc{\operatorname{fc}}
\def\Vsphi{{}^{s\Phi}\Nu}
\def\tx{\tilde x}
\def\sphihalfdens{{}^{s\Phi}\Omega^{1/2}}
\def\mf{\operatorname{mf}}
\def\Diff{\operatorname{Diff}}
\renewcommand\sc{\operatorname{sc}}
\newcommand\Diag{\operatorname{Diag}}
\setcounter{secnumdepth}{3}
\newtheorem{lemma}{Lemma}[section]
\newtheorem{prop}[lemma]{Proposition}
\newtheorem{thm}[lemma]{Theorem}
\newtheorem{cor}[lemma]{Corollary}
\newtheorem{result}[lemma]{Result}
\newtheorem*{thm*}{Theorem}
\numberwithin{equation}{section}
\theoremstyle{remark}
\newtheorem{rem}[lemma]{Remark}
\theoremstyle{definition}
\newtheorem{Def}[lemma]{Definition}
\def\signature#1#2{\par\noindent#1\dotfill\null\\*
{\raggedleft #2\par}}

\renewcommand\div{\operatorname{div}}
\newcommand\Cal{\mathcal}
\newcommand\Lap{\varDelta}
\newcommand\del{\partial}
\renewcommand\a{\alpha}
\renewcommand\b{\beta}
\newcommand\ga{\gamma}
\renewcommand\d{\delta}
\renewcommand\th{\theta}
\newcommand\e{\epsilon}
\newcommand\s{\sigma}
\renewcommand\o{\omega}
\newcommand\Sig{\Sigma}
\newcommand\CI{{\Cal C}^{\infty}}
\newcommand\RR{{\Bbb R}}
\newcommand\ZZ{{\Bbb Z}}
\newcommand\CC{{\Bbb C}}
\newcommand\xa{x_{\alpha}}
\newcommand\ab{\alpha \beta}
\newcommand\Ma{M_{\alpha}}
\newcommand\Va{V_{\alpha}}
\newcommand\Mb{M_{\beta}}
\newcommand\Hab{H_{\alpha\beta}}
\newcommand\bT{{}^bT}
\newcommand\bcT{{}^bT^*}
\newcommand\bw{\bigwedge}
\newcommand\bbw[1]{{}^b\bw{}^{#1}}
\newcommand\Xe{X_{\epsilon}}
\newcommand\ba{\beta \alpha}
\newcommand\ha{\frac12}
\newcommand\ra{\rho_{\alpha}}
\newcommand\be{{}^b\eta}
\newcommand\gab{\ga_{\ab}}
\newcommand\Spin{\text{Spin\,}}
\newcommand\ind{\text{ind\,}}
\newcommand\Sdet{\text{Sdet\,}}
\newcommand\se{s_{\e}}
\newcommand\gep{g_{\e}}
\newcommand\bg{\beta \gamma}
\newcommand\ag{\alpha \gamma}

\newcommand\etag{\eth_{\Cal G}}
\newcommand\tr{\operatorname{tr}}
\newcommand\proj{\operatorname{proj}}
\newcommand\gi{\operatorname{generalized inverse of}}
\newcommand\olS{\overline S}

\newcommand\bH[1]{{}^b\kern-1pt H^{#1}}
\newcommand\bL{{}^b\kern-1pt L}
\newcommand\rest{\restriction}
\newcommand\nul{\operatorname{null}}
\newcommand\clos{\operatorname{clos}}
\newcommand\ins[1]{\overset\circ{#1}}
\newif\ifwantold
\newcommand\old[1]{\ifwantold #1\else\relax\fi}
\newif\ifhide
\let\pa\del
\newcommand\sola{{\frak o}}
\newcommand\SO{\operatorname{SO}}
\newcommand\Xo{X_{\epsilon_0}}
\newcommand\Hba{H_{\beta\alpha}}
\newcommand\APS{\operatorname{APS}}
\newcommand\AS{\operatorname{AS}}
\newcommand\Ch{\operatorname{Ch}}
\newcommand\Refsect[1]{\S\ref{#1}}
\newcommand\Mwhere{\text{ where }}
\newcommand\Mnear{\text{ near }}
\newcommand\lang{\big\langle}
\newcommand\rang{\big\rangle}
\newcommand\lb{\operatorname{lb}}
\newcommand\rb{\operatorname{rb}}
\newcommand\bF{\operatorname{bf}}
\newcommand\Ls{L^{\#}}
\def\ilabel#1{\label{#1}}

\def\thebibliography#1{\section*{References}\list
 {[\arabic{enumi}]}{\settowidth\labelwidth{[#1]}\leftmargin\labelwidth
 \advance\leftmargin\labelsep
 \usecounter{enumi}}
 \def\newblock{\hskip .11em plus .33em minus .07em}
 \sloppy\clubpenalty4000\widowpenalty4000
 \sfcode`\.=1000\relax}
\let\endthebibliography=\endlist

\def\onto{\to}
\def\hz{\hat z}
\def\hzz{\hat z_0}
\def\hzp{\hat z'}
\def\pa{\partial}
\def\({\left(}
\def\){\right)}
\def\hV{\hat V_i *}
\def\ol{\overline}
\def\zetahat{\hat \zeta'}
\def\modzeta{|\zeta'|}
\def\zhat{\hat z}
\def\Sch{Schr\"odinger }
\def\bRN{\overline {{\Bbb R}^N}}
\def\scTX{{}^{sc} T^*_{\pa X} X}
\def\scRN{{}^{sc} T^*_{S^{N-1}} \overline {{\Bbb R}^N}}
\def\tscRN{{}^{sc} T^* \overline {{\Bbb R}^N}}
\def\scpsdo#1#2{{}^{sc}\Psi^{#1, #2}_{cl}(\overline {{\Bbb R}^N})}
\def\scbs#1#2{\sigma_{#2}}
\def\scwf{{}^{sc} \text{WF}}
\def\legsym#1{{}^{sc} \sigma^{#1} }
\def\half{{\frac1{2}}}
\def\modS{\text{ mod } {\Cal S} }
\def\olV{ {\overline{V}}}
\def\hypp{ \hat y''}
\def\zetapy{ {\zeta''}^{\perp}}
\def\Id{\operatorname{Id}}
\def\bgamma{\overline \gamma}
\def\bL{\overline L}
\def\Gs{{G^\sharp}}
\def\Gsl{G^\sharp(\lambda)}
\def\Gt{\tilde G}
\def\Gh{\hat G}
\def\exL{ \overline L_{k, \alpha}}
\def\Id{\operatorname{Id}}
\def\RN{{\Bbb R}^N}
\def\lam{\lambda}
\def\Rml{R_-(\lambda)}
\def\uka{u_{l, \alpha}}
\def\eka{e_{l, \alpha}}
\def\tphi{\tilde \phi}
\def\utl{u_{\theta, \lambda}}
\def\e{\epsilon}
\def\vka{v_{l, \alpha}}
\def\WF{\operatorname{WF}}
\def\CR{\operatorname{CR}}
\def\wzp{\widehat{z'}}
\def\wzpp{\widehat{z''}}
\def\cl{\operatorname{cl}}
\def\interior{\operatorname{int}}
\def\uol{u_{\omega, \lambda}}
\def\dnu#1{\frac{\partial {#1}}{\partial \nu}}
\def\ac{\operatorname{ac}}
\def\td{\tilde}
\def\ot{\leftarrow}
\def\cdotinfty{\dot C^\infty}
\def\sm{\sigma_{\operatorname{max}}}
\def\lzr{\langle z \rangle}
\def\ep{\epsilon}
\def\Op{\operatorname{Op}}
\def\Rea{\operatorname{Re}}
\def\Imag{\operatorname{Im}}
\def\phg{\operatorname{phg}}
\def\Nu{{\Cal V}}
\def\ang#1{\langle {#1} \rangle}
\def\psido{$\Psi$.d.o. }
\def\psidos{$\Psi$.d.o.s }
\def\resA{(A-\lam)^{-1}}
\def\cl{\operatorname{cl}}
\def\slim{\operatorname{s-lim}}
\def\ord{\operatorname{ord}}
\def\ad{\operatorname{ad}}
\def\Psizcl{\Psi^{\ZZ}_{\cl}}
\def\sc{\operatorname{sc}}
\def\X2sc{X^2_{\sc}}
\def\ilg{\operatorname{ilg}}
\def\inc{\operatorname{in}}
\def\out{\operatorname{out}}
\def\ang#1{\big\langle #1 \big\rangle}
\def\GISwf{ {}^{\operatorname{GIS}} \WF}
\def\op{\operatorname{op}}
\def\Rn{\RR^n}
\def\sigb{\sigma_b}
\long\def\skip#1{}
\def\CIdot{\dot \CI}
\def\symrel#1#2#3{\sigma^{#1}_{\mathcal{#2}, #3}}
\def\sc{\operatorname{sc}}
\def\sf{\operatorname{s\Phi}}
\def\scOh{{}^{\sc}{\Omega^{\half}}}
\def\sfOh{{}^{\sf}{\Omega^{\half}}}
\def\schd{\left|\frac{dx dy}{x^{n+1}}\right|^\half}
\def\xt{\tilde x}
\def\jsc{j_{\sc}}
\newcommand\symbint{\sigma_{\operatorname{int}}}
\newcommand\symbb{\sigma_{\pa}}
\def\Lt{\tilde L}
\def\pit{\tilde \pi}
\def\ybar{\overline{y}}
\def\sfOh{{}^{\sf} \Omega^\half}
\def\Sm{S^{[m]}}
\def\Rp{R_+ (\lam)}
\def\Rm{R_- (\lam)}
\def\Rpm{R_\pm (\lam)}
\def\Rpln#1{R_+^{(#1)}(\lam)}
\def\XXbt{\tilde{\XXb}}
\def\Lh{\hat L}
\def\Lhp{\hat L^+}
\def\Lpl{L_+(\lam)}
\def\Lplh{\hat L_+(\lam)}
\def\Lsl{\Ls(\lam)}
\def\bff{\operatorname{bf}}
\def\kid{K_{\Id}}
\def\Rt{\tilde R}

\title[The resolvent for Laplace-type operators] 
{The resolvent for Laplace-type operators on asymptotically conic spaces}
\author{Andrew Hassell}
\address{Centre for Mathematics and its Applications, Australian National
  University, Canberra ACT 0200 Australia}
\email{hassell@maths.anu.edu.au}
\author[Andr\'as Vasy]{Andr\'as Vasy}
\address{Department of Mathematics, University of California, Berkeley,
  CA 94720, U.S.A.}
\email{andras@math.berkeley.edu}


\begin{abstract}Let $X$ be a compact manifold with boundary, and $g$ a
scattering metric on $X$, which may be either of short range or `gravitational'
long range type. Thus, $g$ gives $X$ the geometric structure of a complete
manifold with an asymptotically conic end. 
Let $H$ be an operator of the form $H =
\Lap + P$, where $\Lap$ is the Laplacian with respect to $g$ and
$P$ is a self-adjoint first order scattering differential operator with
coefficients vanishing at $\pa X$ and satisfying a `gravitational'
condition. We define a symbol calculus for Legendre distributions on
manifolds with codimension two corners and use it to give a
direct construction of the resolvent kernel of $H$, $R(\sigma + i0)$, for
$\sigma$ on the positive real axis. 
In this approach, we do not use the limiting
absorption principle at any stage; instead we construct a parametrix  
which solves the resolvent equation up to a compact error term and then use
Fredholm theory to remove the error term. 
\end{abstract}
\maketitle

\section{Introduction} Scattering metrics are a class of Riemannian 
metrics which describe manifolds which
are geometrically complete, and asymptotically conic at infinity. 
We consider manifolds which have 
one end which is diffeomorphic to $Y \times [1, \infty)_r$, where
$Y$ is a closed manifold, and is metrically asymptotic to $dr^2 + r^2 h$,
where $h$ is a Riemannian metric on $Y$, as $r \to \infty$. The precise
definition is given in Definition~\ref{def-sc-met} below. Examples include
the standard metric and the Schwartzschild metric on Euclidean space. 

In this paper we give a direct construction of the
outgoing resolvent kernel, $R(\sigma + i0) = (H - (\sigma + i0))^{-1}$, for
$\sigma$ on the real axis, where $H$ is a perturbation of the Laplacian with
respect to a scattering metric. The incoming resolvent kernel, $R(\sigma -
i0)$, may be obtained by taking the formal adjoint kernel.

The strategy of the proof is to compactify the space to a compact manifold
$X$ and use the scattering calculus of Melrose, as well as the calculus of
Legendre distributions of Melrose-Zworski, extended by us in \cite{HV2}. The
oscillatory behaviour of the resolvent kernel is analyzed in terms of the
`scattering wavefront set' at the boundary. Using propagation of singularity
theorems for the scattering wavefront set leads to an ansatz for the
structure of the resolvent kernel as a sum of a pseudodifferential term and
Legendre distributions of various types. The calculus of Legendre
distributions allows us to construct a rather precise parametrix for the
resolvent in this class, with an error term $E$ which is compact. Using the
parametrix, we show that one can make a finite rank correction to the
parametrix which makes $\Id + E$ invertible, and thus can correct the
parametrix to the exact resolvent. 

As compared to the method of \cite{HV2},
where the authors previously constructed the resolvent, the construction is
direct in two senses. First, we write down rather explicitly a parametrix
for $R(\sigma + i0)$ and then solve away the error using Fredholm
theory. In \cite{HV2}, by contrast, the resolvent was constructed via the
spectral measure, which itself was constructed from the Poisson
operator. Second, we make no use of the limiting absorption principle; that
is, we work directly at the real axis in the spectral variable rather than
taking a limit $R(\sigma + i\ep)$ as $\ep \to 0$. We then prove a posteriori
that the operator constructed is equal to this limit. 

\  

Let us briefly describe the main result here. We consider an operator $H$ of the
form $H = \Lap + P$ acting on half densities, where $P$ is, in the first
place, a short range 
perturbation, that is, a first order self-adjoint differential
operator with coefficients vanishing to
second order at infinity. (Later, we show that there is a simple extension
to metrics and perturbations of `long range gravitational type', which
includes the Newtonian or Coulomb potential and
metrics of interest in general relativity.) We remark that
the Riemannian half-density $|dg|^{1/2}$ trivializes the half-density
bundle, and operators on functions can be regarded as operators on
half-densities via this trivialization.
Given $\lam > 0$, we solve for a kernel $\Rt(\lam)$ on $X^2$ which satisfies
\begin{equation}
(H - \lam^2) \Rt(\lam) = K_{\Id},
\ilabel{eq-Rt}\end{equation}
where $K_{\Id}$ is the kernel of the identity operator on half
densities. More precisely, we 
consider this equation on $\XXb$, which is the space $X^2$ with the corner
blown up. This allows us to use the scattering wavefront set at the `front
face' (the face resulting from blowing up the corner) 
to analyze singularities, which is an absolutely crucial part of the
strategy. The kernel $\Rt(\lam)$ is also required to satisfy a wavefront set
condition at the front face, which is the analogue of the outgoing
Sommerfeld radiation condition. 

We cannot find $\Rt(\lam)$ exactly in one step, so first we look for an
approximation $G(\lam)$ of it.   
The general strategy is to find $G(\lam)$ which solves away the singularities
of the right hand side, $K_{\Id}$, of \eqref{eq-Rt}. Singularities should be
understood both in the sense of interior singularities and oscillations, or
growth, at the boundary, as measured by the scattering wavefront set. 

The first step is to find a pseudodifferential approximation which solves
away the interior singularities of $K_{\Id}$, which is a conormal
distribution supported on the diagonal. This can be done and removes
singularities except at the boundary of the diagonal, where $H - \lam^2$ is
not elliptic (in the sense of the boundary wavefront set). In fact, the
singularities which remain lie on a Legendrian submanifold $N^* \diag_\bl$
at the boundary of $\diag_\bl$ (see \eqref{Ndiag}). Singularities of
$G(\lam)$ can be expected to propagate in a Legendre submanifold $\Lpl$ which
is the bicharacteristic flowout from the intersection of $N^* \diag_\bl$
and the characteristic variety of $H - 
\lam^2$. (The geometry here is precisely that of the fundamental solution of
the wave operator in $\RR^{n+1}$, which is captured by the intersecting
Lagrangian calculus of Melrose-Uhlmann \cite{Rbm-Uhl:Intersecting}.) This
Legendre has conic singularities at another Legendrian, $\Lsl$, 
which is `outgoing'. Thus, in view of the calculus of Legendre distributions
of Melrose-Zworski and the authors, the simplest one could hope for is that
the resolvent on the real axis is the sum of a pseudodifferential term, an
intersecting Legendre 
distribution associated to $(N^* \diag_\bl, \Lpl)$ and a Legendre conic pair
associated to $(\Lpl, \Lsl)$. This is the case:

\begin{thm}\label{main-result} Let $H$ be a short range perturbation of a
short range scattering metric on $X$. Then, for $\lam > 0$, the outgoing
resolvent $R(\lam^2 + i0)$ lies in the class \eqref{seek}, that is, it is the
sum of a scattering pseudodifferential operator of order $-2$, an
intersecting Legendre distribution of order $-1/2$ associated to $(N^*
\diag_\bl, \Lpl)$ and a Legendrian conic pair associated to $(\Lpl, \Lsl)$
of orders $-1/2$ at $\Lpl$, $(n-2)/2$ at $\Lsl$ and $(n-1)/2$ at the left
and right boundaries. 

If $H$ is of long range gravitational type, then the same result holds
except that the Legendre conic pair is multiplied by a complex power of the
left and right boundary defining functions. 
\end{thm} 

This theorem was already proved in our previous work \cite{HV2}, so it is the
method that is of principal interest here. By comparison with \cite{HV2},
the proof is conceptually much shorter; it does not use any results from
\cite{RBMSpec} or \cite{RBMZw}, though it makes substantial use of machinery
from 
\cite{RBMZw}. But the main point we wish to emphasize is that the proof
works directly on the spectrum and 
nowhere uses the limiting absorption principle, a method of attack that we think
will be useful elsewhere in scattering theory. It seems that things which
are easy to prove with this method are difficult with the limiting absorption
principle, and vice versa. For example, it is immediate from our results
that if $f$ is compactly supported in the interior of $X$, then $u = R(\lam^2 +
i0)f$ is such that $x^{-(n-1)/2} e^{-i\lam/x} u \in \CI(X)$, while it is not
so easy to see that the resolvent is a bounded operator from $x^l L^2$ to
$x^{-l} L^2$ for any $l > 1/2$. Using the limiting absorption principle, it
is the second statement that is much easier to derive (following
\cite{GIS} for example). Thus, we hope that
this type of approach will complement other standard methods in scattering
theory. 

\

In the next section, we describe the machinery required, including the
scattering calculus on manifolds with boundary, the scattering-fibred
calculus on manifolds with codimension two corners, and Legendre
distributions in these contexts. The b-double space, which is a blown up
version of the double space $X \times X$ which carries the resolvent kernel,
is also described. The discussion here is rather concise, but there are more
leisurely treatments in \cite{RBMZw} and \cite{HV2}. 

The third section gives a symbol calculus for Legendre distributions on
manifolds with codimension two corners. This is a straightforward
generalization from the codimension one case. 

The fourth section is the heart of the paper, where we construct the
parametrix $G(\lam)$ for the resolvent kernel. Propagation of singularity
theorems show that the simplest space of functions in which one could hope
to find the resolvent kernel is given by \eqref{seek}. We can in fact
construct a parametrix for the resolvent in this class. In the fifth section
this is extended to long range metrics and perturbations. 

In the final section we show that one can modify the parametrix so that the
error term $E(\lam)$ is such that $\Id + E(\lam)$ is invertible. This is done by
showing that the range of $H - \lam^2$ on the sum of $\CIdot(X)$ and $G(\lam)
\CIdot(X)$ is dense on suitable weighted Sobolev spaces. Thus the parametrix
may be corrected to an exact solution of \eqref{eq-Rt}. Such a result also
shows the absence of positive eigenvalues for $H$. Finally, we show that the
kernel so constructed has an analytic continuation to the upper half plane
and agrees with the resolvent there. 

{\it Notation and conventions. } On a compact manifold with boundary, $X$, we
use $\dCinf(X)$ to denote the class of smooth functions, all of whose
derivatives 
vanish at the boundary, with the usual Fr\'echet topology, and $\dist(X)$ to
denote its 
topological dual. On the radial compactification of $\RR^n$ these correspond
to the space of Schwartz functions and tempered distributions, respectively. The
Laplacian $\Delta$ is taken to be positive. The space $L^2(X)$ is taken with
respect to the Riemannian density induced by the scattering metric $g$. This
density has the form $a \, dx dy/x^{n+1}$ near the boundary, where $a$ is
smooth.

{\it Acknowledgements. } We wish to thank Richard Melrose and Rafe Mazzeo for
suggesting the problem, and for many very
helpful conversations. A.\ H.\ is grateful to the Australian Research
Council for financial support. A.\ V.\ thanks the NSF for partial support,
NSF grant \#DMS-99-70607.

\section{Preliminaries}

\subsection{Scattering calculus}

Let $X$ be a manifold with boundary $\pa X = Y$. Near the boundary we will
write local coordinates in the form $(x,y)$ where $x$ is a boundary defining
function and $y$ are coordinates on $Y$ extended to a collar neighbourhood
of $\pa X$. 

We begin by giving the definition of a scattering metric.  
The precise requirements for the metric
(and many other things besides) are easiest to formulate in terms of a
compactification of the space. Taking the function $x = r^{-1}$ as a
boundary defining function and adding a copy of $Y$ at $x=0$ yields a
compact manifold, $X$, with boundary $\pa X = Y$. Then the definition of
scattering metric is given in terms of $X$ in Definition~\ref{def-sc-met}
below. Regularity statements 
for the metric coefficients are in terms of the 
$\CI$ structure on $X$; this is a strong requirement, being equivalent to
the existence of a complete asymptotic expansion, together
with all derivatives, in inverse powers of $r$ as $r \to \infty$. The
benefit of such a strong requirement is that we get complete asymptotic
expansions for the resolvent kernel, and correspondingly, mapping properties
of the resolvent on spaces of functions with complete asymptotic
expansions. 

\begin{Def}\ilabel{def-sc-met} 
A (short range) scattering metric on $X$ is a Riemannian metric
$g$ in the interior of $X$ which takes the form
\begin{equation}
g = \frac{dx^2}{x^4} + \frac{h'}{x^2},
\label{sc-metric-1}\end{equation}
where $h'$ is a smooth symmetric 2-cotensor on $X$ which restricts to the
boundary to be a metric $h$ on $Y$ \cite{RBMSpec}. A long range scattering
metric is a 
metric in the interior of $X$ which takes the form
\begin{equation}
g = a_{00} \frac{dx^2}{x^4} + \frac{h'}{x^2},
\label{sc-metric-2}\end{equation}
where $a_{00}$ is smooth on $X$, $a_{00} = 1 + O(x)$, and $h'$
is as above \cite{Vasy:Geometric}. If $a_{00} = 1 - cx + O(x^2)$ for some
constant $c$ we call $g$ a gravitational long range scattering metric. 
\end{Def}

{\bf Examples.} Flat Euclidean space has a metric which in polar coordinates
takes the form
$$
dr^2 + r^2 d\omega^2,
$$
where $d\omega^2$ is the standard metric on $S^{n-1}$. Compactifying Euclidean
space as above, we obtain a ball with $x = r^{-1}$ as boundary defining
function, and then the flat metric becomes
$$
\frac{dx^2}{x^4} + \frac{d\omega^2}{x^2},
$$
which is a short range scattering metric. 

The Schwartzschild metric on $\RR^n$ takes the form near infinity
$$
\left( 1 - \frac{2m}{r} \right) dr^2 + r^2 d\omega^2,
$$
which under the same transformation leads to a gravitational long range
scattering metric 
$$
(1 - 2mx) \frac{dx^2}{x^4} + \frac{d\omega^2}{x^2}.
$$
The constant $m = c/2$ is interpreted as the mass in general relativity. 

\

The natural Lie Algebra corresponding to the class of scattering metrics
on $X$ is the scattering Lie Algebra 
$$
\Vsc(X) = \{ V \mid V = x W, \ \text{where } W \text{ is a } \CI \text{ vector
  field on } X \text{ tangent to } \pa X \}.
$$
Clearly this Lie Algebra can be localized to any open set. In the interior
of $X$, it consists of all smooth vector fields, while near the boundary it
is equal to the $\CI(X)$-span of the vector fields $x^2 \pa_x$ and $x
\pa_{y_i}$. Hence it is the space of smooth sections of a vector bundle,
denoted $\Tsc X$, the scattering tangent bundle. Any scattering metric
turns out to be a smooth fibre metric on $\Tsc X$. The dual bundle, denoted
$\sct X$, is called the scattering cotangent bundle; near the boundary,
smooth sections are generated over $\CI(X)$ by $dx/x^2$ and $dy_i/x$. A
general point in $\sct _p X$ can be thought of as the value of a
differential $d(f/x)$ at $p$, where $f \in \CI(X)$, and in terms of local
coordinates $(x,y)$ near $\pa X$ can be written 
$\tau dx/x^2 + \mu_i dy_i/x$,
yielding local coordinates $(x,y,\tau,\mu)$ on $\sct X$ near $\pa X$. 

The scattering
differential operators of order $k$, denoted $\Diffsc ^k(X)$, are those given
by sums of products of at most $k$ scattering vector fields. There are two
symbol maps defined for $P \in \Diffsc^k(X)$. The first is the `usual' symbol map,
denoted here $\symbint^k(P)$, which maps to $S^k(\sct X)/S^{k-1}(\sct X)$,
where $S^k(\sct X)$ denotes the classical symbols of order $k$ on
$\sct X$.
The second is the boundary symbol, $\symbb(P) \in S^k(\sct_{\pa X} X)$, 
which is the {\it full} symbol of $P$ restricted to $x=0$. This is well
defined since the Lie Algebra $\Vsc(X)$ has the property $[\Vsc(X), \Vsc(X)]
\subset x \Vsc(X)$, so commutators of scattering vector fields vanish to an
additional order at the boundary. Dividing the interior symbol
$\symbint^k(P)$ by $|\xi|^k_g$, where $|\cdot|_g$ is the metric on $\sct X$
determined by the scattering metric, we get a function on the sphere bundle
of $\sct X$. This may be combined with the boundary symbol into a joint symbol,
$\jsc^k(P)$, a function on $\Csc(X)$ which is the topological space obtained
by gluing together the sphere bundle of $\sct X$ with the the fibrewise
radial compactification of $\sct_{\pa X} X$ along their common boundary. 

\begin{figure}\centering
\epsfig{file=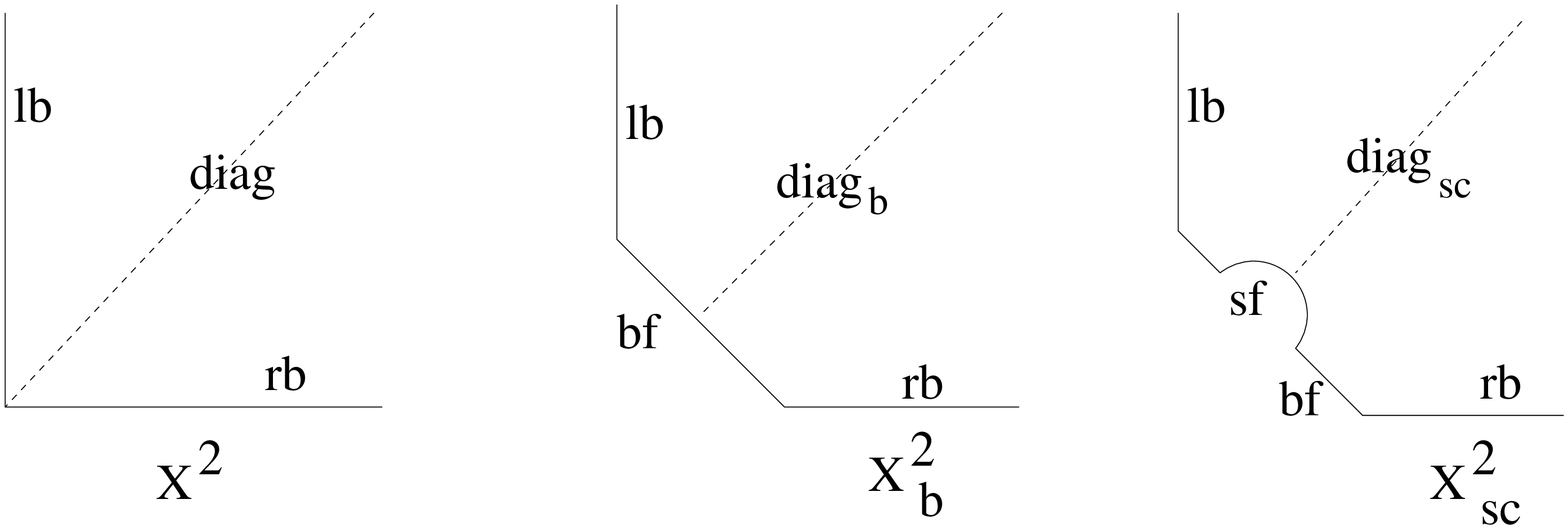,width=12cm,height=4cm}
\caption{The b-double space and the scattering scattering space}
\label{fig:XXb}    
\end{figure}

The scattering pseudodifferential operators are defined in terms of the
behaviour of their Schwartz kernels on the scattering double space 
$\XXsc$, a blown up version of the
double space $X^2$. This is defined by
\begin{equation}\begin{gathered}
\XXb = [X^2; (\pa X)^2] \qquad \text{and} \\
\XXsc = [\XXb; \pa \diag_\bl ],
\end{gathered}\end{equation}
and $\diag_\bl$ is the lift of the diagonal to $\XXb$. The lift of
$\diag_\bl$ to $\XXsc$ is denoted $\diag_\scl$. The blowup notation $[; ]$
is that of Melrose: see \cite{RBMCalcCon} or \cite{tapsit}. The boundary
  hypersurfaces are 
labelled lb, rb, bf and sf; see figure~\ref{fig:XXb}. The scattering
pseudodifferential operators of order $k$, acting on half densities,
$\Psisc^k(X)$, are those given by $\KDsc$-valued distribution kernels which
are classical conormal at $\diag_\scl$, of order $k$, uniformly to the
boundary, and rapidly vanishing at lb, rb, bf. (Here $\KDsc$ is the pullback
of the bundle $\pi_l^* \scOh(X) \otimes \pi_r^* \scOh(X)$ over $X^2$ to $\XXb$.)
The space $\Psisc ^{k,l}(X)$ is defined 
to be $x^l \Psisc ^{k}(X)$. 

The joint symbol map $\jsc^k$ extends from $\Diffsc^k(X)$ to $\Psisc^k(X)$
multiplicatively,
$$
\jsc^k(A) \cdot \jsc^m(B) = \jsc^{k+m}(AB),
$$
and such that there is an exact sequence
$$
\begin{CD}
0 @>>> \Psisc^{m-1,1}(X) @>>> \Psisc^{m}(X) @>{\jsc^{m}}>> 
\CI(\Csc(X)) @>>> 0.
\end{CD}
$$
A scattering pseudodifferential operator $A \in \Psisc^k(X)$ is said to be
elliptic at a point $q \in \Csc(X)$ if $\jsc^k(A)(q) \neq 0$. It is said to 
have 
elliptic interior symbol (boundary symbol) if $\jsc^k(A)$ does not vanish at
fibre-infinity (spatial infinity), and is said to be totally elliptic if
$\jsc^k(A)$ vanishes nowhere. The characteristic variety of $A$,
$\Sigma(A)$, is the zero set of $\jsc^k(A)$. 

The scattering wavefront set of a tempered distribution $u \in \dist(X)$
(the dual space of $\CIdot(X)$, the space of smooth functions on $X$
vanishing with all derivatives at the boundary) is the closed subset of
$\Csc(X)$ whose complement is
\begin{equation}
\scwf(u)^\complement = \{ q \in \Csc(X) \mid \ \exists A \in \Psisc^0(X)
\text{ elliptic at } q \text{ such that } Au \in \CIdot(X) \}.
\end{equation}
The interior part of the wavefront set (at fibre-infinity) is a familiar
object, the standard wavefront set introduced by H\"ormander (except that
each ray of the standard wavefront set is thought of here as a point in the
cosphere bundle). 
In this paper we are mostly interested in the part of the scattering
wavefront set at spatial infinity. In fact, the operators $H$ we shall study
will have elliptic interior symbol, uniformly to the boundary, 
so in view of the next theorem, solutions
of $(H - \sigma) u = 0$ must have wavefront set contained in the
part of $\Csc(X)$ at spatial infinity, which we denote $K$ 
(that is, $K = \sct_{\pa X} X$).

There is a natural contact structure on $K$ induced by the symplectic form
$\omega$ on $T^* X$. Writing $\omega$ in terms of the rescaled cotangent
variables $\tau, \mu$ and contracting with the vector $x^2 \pa_x$ yields the
1-form
$$
\chi = \iota_{x^2 \pa_x} \omega = d\tau + \mu \cdot dy,
$$
which is nondegenerate, and therefore a contact form. A change of boundary
defining function $x' = a x$ changes $\chi$ by a factor $a^{-1}$, so the
contact {\it structure} is totally well-defined. Given a Hamiltonian $h$ on
$K$, the Hamiltonian vector field on $K$ determined by the contact form is
$$
\frac{\pa h}{\pa \mu_i} \frac{\pa}{\pa y_i} + \Big(
-\frac{\pa h}{\pa y_i} + \mu_i \frac{\pa h}{\pa \tau} \Big) \frac{\pa }{\pa
  \mu_i} + \Big( h - \mu_i \frac{\pa h}{\pa \mu_i} \Big) \frac{\pa}{\pa
  \tau}.
$$
This is the same as $x^{-1} V_h$ restricted to $x=0$, where $V_h$ is the
Hamilton vector field on $\sct X$ induced by $h$. 
Integral curves of this vector field are called bicharacteristics of $h$ (or
of $A$, if $h$ is the boundary symbol of $A$). 

Under a coordinate change, the variables $\tau$ and $\mu$ change according
to
$$
\tau' = a \tau, \quad \mu' = a \mu \frac{ \pa y}{\pa y'} - \tau \frac{\pa
  a}{\pa y'} \qquad x' = a x.
$$
Since $a > 0$, this shows that the subset 
\begin{equation}
K_- = \{ (y, \tau, \mu) \in K \mid \tau \leq 0 \}
\ilabel{Kminus}
\end{equation} 
is invariantly defined. This is important in the definition of the outgoing
resolvent, see \eqref{outgoing}. 

The boundary part of the scattering wavefront set behaves very much as the
interior wavefront set part behaves, and in particular we have a propagation
of singularities result for operators of real principal type:

\begin{thm}[Melrose] Suppose $A \in \Psisc^{k}(X)$ has elliptic interior
symbol, and 
real boundary symbol. Then for
$u \in \Cinf(X^o) \cap \dist(X)$,
we have 
\begin{equation}
\WFsc(Au) \subset \WFsc(u), 
\ilabel{prop-1}
\end{equation}
\begin{equation}
\WFsc(u) \setminus \WFsc(Au) \subset \Sigma(A),  
\ilabel{prop-2}\end{equation} 
and 
\begin{equation}\begin{gathered}
\WFsc(u) \setminus \WFsc(Au) 
\text{ is a union of } \\ \text{ maximally
extended bicharacteristics of } A \text{ inside } 
\Sigma(A) \setminus \WFsc(Au).
\end{gathered}\ilabel{prop-3}\end{equation}
\end{thm}

Thus, if $Au = 0$, then $\scwf(u) \subset K$ and consists of a union of
maximally extended bicharacteristics of $A$ inside $\Sigma(A)$. 

As well as a boundary principal symbol defined on $K$, scattering
pseudodifferential operators also have a boundary subprincipal symbol. This
is the $O(x)$ term of the full symbol at the boundary when the operator is
written in Weyl form. It is important to keep in mind that it depends on a 
choice of product structure at the boundary; it does not enjoy quite the
same invariance properties as does the standard (interior) subprincipal
symbol.  A practical formula to use for differential 
operators with symbol in left-reduced form, ie, such that 
$$
\sigma \left( a(x,y) (x^2 D_x)^j (x D_y)^\alpha \right) = a(x,y) \tau^j
\mu^\alpha, \quad D = -i \pa, \quad \alpha =
(\alpha_1, \dots, \alpha_{n-1}),
$$
is that for $\sigma(P) = p(y, \tau, \mu) + x q(y, \tau, \mu) + O(x^2)$,
the boundary subprincipal symbol of $P$ is given by
\begin{equation}
\sigma_{\sub}(P) = q + \frac{i}{2} \left( \frac{\pa^2 p}{\pa y_i \pa \mu_i}
- (n-1) \frac{ \pa p}{\pa \tau} + \mu_i \frac{ \pa^2 p}{\pa \mu_i \pa \tau}
\right).
\ilabel{subpr-symb}\end{equation}

\begin{lemma}\label{lem:sub-pr} 
Let $g$ be a short range scattering metric, let $x$ be a
boundary defining function with respect to which 
$g = dx^2/x^4 + h'/x^2$, and let $H$ be a short range perturbation of the
Laplacian with respect to $g$. Then in local coordinates $(x,y)$, the 
sub-principal symbol of $H$ vanishes at $\mu = 0$. 
\end{lemma}

\begin{proof} The operator $H$ may be written
$$
H = (x^2 D_x)^2 + (n-1) i x^3 D_x + x^2 \Lap_\pa + a_{ij}(x,y) x^3 D_{y_i}
D_{y_j} + Q, \quad Q \in x^2 \Diffsc(X). 
$$
Thus the left-reduced symbol as above is
$$
\sigma(H) = \tau^2 + h_{ij}\mu_i \mu_j + x\left( (n-1) i \tau + a_{ij} \mu_i
\mu_j \right) + O(x^2).
$$
Hence the sub-principal symbol is 
\begin{equation}
\sigma_{\sub}(H) = i \frac{\pa h_{ij}}{\pa y_i} \mu_j + a_{ij} \mu_i \mu_j,
\ilabel{mu=0}\end{equation}
which vanishes when $\mu = 0$. 
\end{proof}

We now define the gravitational condition for the perturbation $P$. 

\begin{Def} A first order scattering differential operator $P$ on $X$ is
said to be 
short range if it lies in $x^2 \Diffsc^1(X)$, and long range if it lies in $x
\Diffsc^1(X)$. Let $g$ be a scattering metric and $x$ a boundary defining
function with respect to which $g$ takes the form \eqref{sc-metric-1} or
\eqref{sc-metric-2}. $P$ is said to be 
of long range gravitational type with respect to $g$ if it has the form
$$
P = x \left( \sum_i a_i (x \pa_{y_i}) + b x^2 \pa_x + c \right),
$$
near $x=0$, where $a_i$, $b$ and $c$ are in $\CI(X)$, and for some
constants $b_0$ and $c_0$, $b = b_0 + O(x)$ and $c = c_0  + O(x)$. 
\end{Def}

The point of the short range condition is that then the subprincipal symbol of
both $H = \Lap +  P$ vanishes at the radial sets $\mu = 0, \tau = \pm \lam$
of $H - \lam^2$, whilst in the long range gravitational case, the subprincipal
symbol is constant. In the general long range case, the subprincipal symbol
is an arbitrary function on the radial set, which causes some inconvenience
(but not insuperable difficulties) in constructing the parametrix for
$(H-\lam^2 - i0)^{-1}$. 

\subsection{Legendre distributions} An important special case that occurs
often is that $\scwf(u)$ is a Legendre submanifold, or union of Legendre
submanifolds, of $K$; moreover, in many cases, $u$ is a Legendre
distribution, which means that it has a WKB-type expansion, the product of a
oscillatory and smooth term, as discussed below, which makes it
particularly amenable to analysis. 

We let $\dim X = n$, so that $\dim K = 2n-1$. 
A Legendre submanifold of $K$ is a submanifold $G$ of dimension $n-1$ such that $\chi
\restriction G = 0$. Such submanifolds have several nice properties. One is
that if a Hamiltonian, $h$, is constant on $G$ then its Hamilton vector
field is tangent to $G$. Another is that Legendre submanifolds may be
generated in the following way: If $F$ is a submanifold of dimension $n-2$,
such that $\chi$ vanishes on $F$, and if the Hamilton vector field of $h$ is
nowhere tangent to $F$, then the union of bicharacteristics of $h$ passing
through $F$ is
(locally) a Legendre submanifold. 

Let $G$ be a Legendre submanifold, and let 
$q \in G$. A local (nondegenerate) parametrization of $G$ near $q$ is a
function $\phi(y,v)$ defined in a neighbourhood of $y_0 \in Y$ and $v_0 \in
\RR^k$, such that $d_v \phi = 0$ at $q' = (y_0, v_0)$, $q = (y,
d_{(x,y)}(\phi/x))$ at $q'$, 
$\phi$ satisfies the nondegeneracy hypothesis
$$
d\Big( \frac{ \pa \phi}{\pa v_i} \Big) \text{ are linearly independent at }
C = \{ (y,v) \mid d_v \phi = 0 \}, 1 \leq i \leq k,
$$
and near $q$,
\begin{equation}
G = \{ \big(y, d_{(x,y)} \left( \frac{\phi}{x} \right) \big) \mid (y,v) \in C \}.
\ilabel{correspondence}\end{equation}
A Legendre distribution of order $m$ associated to $G$ is a half-density of
the form $u = (u_0 + \sum_{j=1}^N u_j) \nu$, where $\nu$ is a smooth section
of the scattering half density bundle, $u_0 \in \CIdot(X)$, and
$u_j$ is supported in a coordinate patch $(x,y)$ near the boundary, with an
expression 
$$
u_j = x^{m+n/4-k/2} \int_{\RR^k} e^{i\phi_j(y,v)/x} a_j(x,y,v) \, dv,
$$
where $\phi_j$ locally parametrizes $G$ and $a_j \in \CI(X \times \RR^k)$,
with compact support in $v$. Melrose and Zworski showed that $u_j$ can be
written with respect to any local parametrization, up to an error in
$\CIdot(X)$. The set of such half-densities is denoted $I^m(X, G;
\scOh)$. The scattering wavefront set of $u \in I^m(X, G; \scOh)$ is
contained in $G$.  

An intersecting Legendre distribution is associated to a pair of Legendre
submanifolds, $\Lt = (L_0, L_1)$, where $L_1$ is a manifold with boundary such
that $L_0$ and $L_1$ intersect cleanly at $\pa L_1$. A local parametrization
of $(L_0, L_1)$ near $q \in L_0 \cap L_1$ is a function $\phi(y,v,s)$
defined in a neighbourhood of $q' = (y_0,v_0,0)$ in $Y \times \RR^k \times [0,
\infty)$ such that $d_v \phi = 0$ at $q'$, $q = (y, d_{(x,y)}(\phi/x))(q')$,
$\phi$ satisfies the nondegeneracy hypothesis
$$
ds, \ d\phi, \ \text{ and }
d\Big( \frac{ \pa \phi}{\pa v_i} \Big) \text{ are linearly independent at }
q', \ 1 \leq i \leq k,
$$
and near $q$,
$$ \begin{gathered}
L_0 = \{ \big(y, d_{(x,y)} \left( \frac{\phi}{x} \right) \big) \mid s=0, d_v
\phi = 0 \}, \\ 
L_1 = \{ \big(y, d_{(x,y)} \left( \frac{\phi}{x} \right) \big) \mid s\geq 0,
d_s \phi = 0, \ d_v \phi = 0 \}. 
\end{gathered}
$$
A Legendre distribution of order $m$ associated to $\Lt$ is a 
half-density of
the form $u = u_0 + (\sum_{j=1}^N u_j)\nu$, where $\nu$ is a smooth
scattering half-density, $u_0 \in I^m_c(X, L_1; \scOh) +
I^{m+1/2}(X, L_0; \scOh)$ (the subscript $c$ indicates that the microlocal
support does not meet the boundary of $L_1$), and
$u_j$ is supported in a coordinate patch $(x,y)$ near the boundary, with an
expression 
$$
u_j = x^{m+n/4-(k+1)/2} \int_0^\infty ds \int e^{i\phi_j(y,v,s)/x} a_j(x,y,v,s) 
\, dv  ,
$$
where $\phi_j$ locally parametrizes $(L_0, L_1)$ and 
$a_j \in \CI(X \times \RR^k \times [0, \infty))$,
with compact support in $v$ and $s$. Again, $u_j$ can be
written with respect to any local parametrization, up to an error in
$\CIdot(X)$. The set of such half-densities is denoted 
$I^m(X, \Lt; \scOh)$. The scattering wavefront set of $u \in I^m(X, \Lt;
\scOh)$ is contained in $L_0 
\cup L_1$.

\

A Legendre distribution associated to a conic Legendrian pair is associated
to a pair of Legendre submanifolds $\Gt = (G, \Gs)$ where $\Gs$ is a
projectable Legendrian (that is, the projection from $\sct X$ to $Y$ is a
diffeomorphism restricted to $\Gs$) 
and $G$ is an open Legendrian submanifold such that
$\overline{G} \setminus G$ is contained in $\Gs$ and $\overline{G}$ has at
most a conic singularity at $\Gs$. We further assume that $\tau \neq 0$ on
$\Gs$, so that we may change coordinates to a new boundary defining function
such that $\Gs$ is parametrized by the phase function $1$. In these
coordinates, the
condition that $\overline{G}$ has a conic singularity at $\Gs$ means that 
$\overline{G}$ lifts to a smooth
submanifold with boundary, $\Gh$, on the blown-up space
\begin{equation}
[\sct X; \ \{ x = 0, \mu = 0 \}],
\ilabel{Ghat}\end{equation}
intersecting the front face of \eqref{Ghat} transversally. In local
coordinates $(x,y,\tau,\mu)$, coordinates near the front face are 
\begin{equation}
x/|\mu|, \ y, \ \tau, \ |\mu| \text{ and } \hat \mu,
\ilabel{hat-coords}\end{equation}
and we require that $\Gh$ is given by the vanishing of
$n$ smooth functions of these variables with linearly independent
differentials, and that $d|\mu| \neq 0$ at $\pa 
\Gh$.

A local parametrization
of $\Gt$ near $q \in \overline{G} \cap \Gs$ is a function $\phi(y,v,s) =
1 + s \psi(y,v,s)$
defined in a neighbourhood of $q' = (y_0,v_0,0)$ in $Y \times \RR^k \times [0,
\infty)$ such that $\phi_0$ parametrizes $\Gs$ near $q$, 
$d_v \phi = 0$ at $q'$, $q = (y, d_{(x,y)}(\phi/x))(q')$, $\phi$
satisfies the nondegeneracy hypothesis
$$
ds, \ d\psi, \ \text{ and }
d\Big( \frac{ \pa \psi}{\pa v_i} \Big) \text{ are linearly independent at }
q', \ 1 \leq i \leq k,
$$
and near $q$,
$$ 
\Gh = \{ (0, y, -\phi, s d_y \psi, \widehat{d_y \psi}) 
\mid d_v \phi = 0, \ d_s \psi = 0 \},
$$
in the coordinates \eqref{hat-coords}. 
A Legendre distribution of order $(m,p)$ associated to $(G, \Gs)$ is a 
half-density of
the form $u = u_0 + (\sum_{i=1}^N u_i)\nu$, where $\nu$ is as above, 
$u_0 \in I^m_c(X, G; \scOh) +
I^{p}(X, \Gs; \scOh)$ (the subscript $c$ indicates that the microlocal
support does not meet $\Gs$), and
$u_j$ is supported in a coordinate patch $(x,y)$ near the boundary, with an
expression 
$$
u_j =  \int_0^\infty ds \int e^{i\phi_j(y,v,s)/x} a_j(y,v,x/s,s) \Big
( \frac{x}{s} \Big)^{m+n/4-(k+1)/2} s^{p + n/4 - 1} 
\, dv  ,
$$
where $\phi_j$ locally parametrizes $(G, \Gs)$ and 
$a_j \in \CI(X \times \RR^k \times [0, \infty) \times [0, \infty))$,
with compact support in $v$, $x/s$ and $s$. Here $u_j$ can be
written with respect to any local parametrization, up to an error in
$I^p(X, \Gs; \scOh)$. The set of such half-densities is denoted 
$I^{m,p}(X, \Gt; \scOh)$. 
The wavefront set of $u \in I^m(X, G; \scOh)$ is contained in $G
\cup \Gs$.

\subsection{Codimension 2 corners} In this subsection we briefly review the
extension of the theory of Legendre distributions to manifolds with
codimension 2 corners and fibred boundaries given in \cite{HV2}. 

Let $M$ be a compact manifold with codimension 2 corners. The boundary
hypersurfaces will be labelled $\mf, H_1, \dots, H_d$, where the $H_i$ are
endowed with fibrations $\pi_i : H_i \to Z_i$ to certain closed manifolds
$Z_i$ and mf (the `main face') is given the trivial fibration
$\operatorname{id} : \mf \to \mf$. The collection of fibrations is denoted
$\Phi$. It is assumed that $H_i \cap H_j = \emptyset$ if $i \neq j$.
It is also assumed that the fibres of $\pi_i$ intersect $H_i \cap \mf$
transversally  and therefore induce a fibration from $H_i \cap \mf \to
Z_i$. Further, it is assumed that a total boundary defining function $x$ is
given, which is distinguished up to multiplication by positive functions
which are constant on the fibres of $\pa M$. 

Near $H \cap \mf$, where $H = H_i$ for some $i$, there are coordinates $x_1,
x_2, y_1, y_2$ such that $x_1$ is a boundary defining function for $H$,
$x_2$ is a boundary defining function for mf, $x_1x_2 = x$, and the
fibration on $H$ takes the form
$$
(y_1, x_2, y_2) \mapsto y_1.
$$

Associated with this structure is a Lie Algebra of vector fields
$$
\Vsf(M) = \{ V \mid V \in \CI, V \text{ is tangent to } \Phi \text{ at } \pa M,
V(x) = O(x^2) \}.
$$
This is the space of smooth sections of a vector bundle, denoted $\Tsf M$. 
The dual space is denoted $\sfT M$. A point in $\sfT _p M$ may be thought of
as a differential $d(f/x)$ at $p$, where $f$ is a smooth function on $M$
constant on the fibres at $\pa M$. A basis for $\sfT _p M$, for $p \in M$
near $\mf \cap H$, is given by 
$dx/x^2$, $dx_1/x$, $dy_1/x$, $dy_2/x^2$. Writing $q \in \sfT M$ as
$$
q = \tau \frac{dx}{x^2} + \tau_1 \frac{dx_1}{x} + \mu_1 \cdot
\frac{dy_1}{x} + \mu_2 \cdot \frac{dy_2}{x_2} 
$$
gives coordinates 
\begin{equation}
(x_1, \ x_2, \ y_1, \ y_2, \ \tau, \ \tau_1, \ \mu_1, \ \mu_2)
\ilabel{sfT-coords}\end{equation} 
on $\sfT M$ near $\mf \cap H$. 

The differential operators of order at most $k$ generated over $\CI(M)$ by
$\Vsf(M)$ are denoted $\Diffsf ^k(M)$. Near the interior of mf, the Lie
Algebra $\Vsf(M)$ localizes to the scattering Lie Algebra $\Vsc(\Mt)$, where
$\Mt$ denotes the noncompact manifold with boundary $M \setminus \cup_i
H_i$. Consequently, we have a boundary symbol $\symbb(P)$, $P \in \Diffsf
^k(M)$ taking values in $S^k(\sct _{\pa \Mt} \Mt)$ over the interior of
mf. In fact the symbol extends to an element of $S^k (\sfT _{\mf} M)$
continuous up to the boundary of mf.

For each fibre $F$ of $H$, there is a subbundle of $\Tsf M$ consisting of
all vector fields vanishing at $F$. The annihilator subbundle of $\sfT M$ is
denoted $\sct (H; F)$ since it is isomorphic to the cotangent space of the
fibre. The quotient bundle, $\sfT M / \sct (H; F)$ is denoted $\sfN Z_i$
since it is the pullback of a bundle over $Z_i$. The fibration $\pi_i$
induces a fibration
\begin{equation}\begin{gathered}
\pit_i : \sfT _{\mf \cap H} M \to \sfN Z_i \\
(y_1, y_2, \tau, \tau_1, \mu_1, \mu_2) \mapsto (\tau, y_1, \mu_1).
\end{gathered}\ilabel{pit}
\end{equation}

We next describe three types of contact structures associated with the
structure of 
$M$. Since $\Vsf (M)$ is locally the scattering structure near the interior
of mf, there is an induced contact structure on $\sfT M$ over the interior
of mf. In local coordinates, the contact form looks like
$$
\chi = \iota_{x^2 \pa_x} (\omega) = d\tau + \tau_1 dx_1 + \mu_1 dy_1 + x_1
\mu_2 dy_2.
$$
We see from this that at $x_1 = 0$, $\chi$ is degenerate. However,
restricted to $\mf \cap 
H$, $\chi$ is the lift of a form $\chi_{Z_i}$ on $\sfN (Z_i)$,
namely $d\tau + \mu_1 dy_1$, which is nondegenerate on $\sfN (Z_i)$. This
determines 
our second type of contact structure (one for each $i$). 
The third type of contact structure is that on $\sct
_{\pa F} (H_i; F)$ induced by $\Vsf (M)$ for each fibre $F$ of $H_i$, since
it restricts to the scattering vector fields on each fibre. In local
coodinates, this looks like $d\tau_1 + \mu_2 dy_2$. 

Using these three contact structures we define Legendre submanifolds and
Legendre distributions.

\begin{Def}\label{M-Leg-submfld}
A  Legendre submanifold $G$ of $\sfT M$ is a Legendre
submanifold of $\sfT_{\mf}M$ which is transversal to $\sfT_{H_i \cap
\mf}M$ for each $H_i$,
for which the map \eqref{pit} induces a fibration from $\pa G$ to $G_1$,
where $G_1$ is a Legendre
submanifold of $\sfN Z_i$, whose fibers are Legendre submanifolds
of $\sct_{\partial F} F$.
\end{Def}

A projectable Legendrian (one such that the projection from $\sfT _{\mf} M
\to \mf$ is a diffeomorphism when restricted to $G$) is always of the form 
$$
\text{ graph } \big( d \big(\frac{\phi}{x} \big) \big) = \{ \big(\ybar,
d\big(\frac{\phi(\ybar)}{x}\big)\big)  \mid \ybar \in \mf \}
$$
for some smooth function $\phi$ constant on the fibres of $\pa M$. We then
say that $\phi$ parametrizes $G$. In
general, let $G$ be a Legendre submanifold of $\sfT M$, and let $q \in
G$. If $q$ lies above the interior of mf, then a local parametrization of
$G$ near $q$ is as described in the previous subsection, so consider
$q \in \pa G$ lying in $\sfT _{\ybar_0} M$, where $\ybar_0 \in \mf \cap H$. 
A local (nondegenerate) parametrization of $G$ near $q$ is a
function $\phi(x_1, y_1, y_2, v, w)$ of the form 
\begin{equation}\label{eq:fib-ph-9}
\phi(x_1, y_1, y_2,v, w)=\phi_1(y_1,v)+x_1\phi_2(x_1,y_1,y_2,v, w),
\end{equation}
defined in a neighbourhood of $q' = (\ybar_0,v_0,w_0) \in \mf \times
\RR^{k_1 + k_2}$, such that $d_{v,w} \phi = 0$ at $q'$, 
$$
q = \big(\ybar_0, d\big(\frac{\phi(\ybar_0)}{x}\big)\big) \text{ at } q'
$$ 
in local coordinates \eqref{sfT-coords}, $\phi$
satisfies the nondegeneracy hypothesis at $q'$
\begin{equation}\label{eq:corner-20}
d_{(y_1,v)}\frac{\pa \phi_1}{\partial {v_j}}, j=1,\ldots, k
\text{ and } d_{(y_2,w)} \frac{\pa \phi_2}
{\partial {w_{j'}}},\ j'=1,\ldots,k' \text{ linearly independent,}
\end{equation}
and near $q$,
\begin{equation}
G = \{ \big(\ybar, d\big(\frac{\phi(\ybar)}{x}\big)\big) \mid 
d_v \phi = d_w \psi = 0 \}.
\end{equation}
A Legendre distribution of order $(m; r_1, \dots r_d)$ associated to $G$ is
a half-density 
such that for any $\upsilon_i \in \CI(M)$ whose support does not intersect
$H_k$, for $k \neq i$, $\upsilon_i u$ is of
the form $u = u_0 + (\sum_{j=1}^N u_j + \sum_{j=1}^M u_j')\nu$, where 
$\nu$ is a smooth section of the half-density bundle $\sfOh$ induced by
$\sfT M$, 
$u_0 \in \CIdot(X; \sfOh)$, and $u_j$, $u'_j$ have expressions 
\begin{equation}\begin{gathered}\label{eq:dist-1}
u_j(x_1, x_2,y_1, y_2)
= \int e^{i\phi_j(x_1, y_1, y_2,v,w)/x} a_j(x_1, x_2, y_1, y_2,v,w)\\
x_2^{m-(k+k')/2+N/4}x_1^{r_i-k/2+N/4-f_i/2}\,dv\,dw 
\end{gathered}\end{equation}
with $N=\dim M$, $a_j\in\Cinf_c([0,\epsilon)\times
U\times\RR^{k+k'})$, $U$ open in $\mf$, $f_i$ the dimension of the fibres
of $H_i$ and 
$\phi_j$ a phase function parametrizing a Legendrian $G$ on
$U$, and 
\begin{equation}
u'_j(x_1, y_1,z)
= \int e^{i\psi_j(y_1,w)/x} a_j(x,y_1,z,w)
x^{r_i-k/2+N/4-f_i/2}\,dw 
\ilabel{eq:dist-2}
\end{equation}
with $N=\dim M$, $a_j\in\Cinf_c([0,\epsilon)\times
U\times\RR^{k})$, $U$ open in $H$, $f_i$ as above,
$\psi_j$ a phase function parametrizing the Legendrian $G_1$. 

\begin{Def}\label{M-leg-submfld-conicpts}
A  Legendre pair with conic points, $(G,\Gs)$,
in $\sfT M$ consists of two  Legendre submanifolds $G$ and $\Gs$
of $\sfT M$ which form an intersecting pair with conic points in $\sct_{\td M} \td M$ 
such that for each $H_i$ the
fibrations of $G$ and $\Gs$ induced by \eqref{pit} have the
same Legendre submanifold $G_1$ of $\sfN Z_i$ as base and for
which the fibres are intersecting pairs of 
Legendre submanifolds with conic points
of $\sct_{\partial F} F$. 
\end{Def}

The Legendrian $\Gs$ is required to be projectable, so it parametrized by a
phase function $\phi(\ybar)$ which is constant on the fibres of $\pa
M$. Thus, $x' = x/\phi$ is another admissible total boundary defining
function. With respect to $x'$, $\Gs$ is
parametrized by the function $1$. Thus, without loss of generality we may
assume that coordinates have been chosen so that $\Gs$ is parametrized by
$1$. This simplifies the coordinate form of the blowup
\eqref{Ghat}. Coordinates near $\pa \Gh$ then are 
\begin{equation}
\ybar = (x_1, y_1, y_2), \ \tau,\ \tau_1/|\mu_2|, \ \mu_1/|\mu_2|, \ 
\muh_2 = \mu_2/|\mu_2|, \text{ and } |\mu_2|, 
\ilabel{Ghat-coords}\end{equation}
the last of which is a boundary defining function for $\Gh$ (see
\cite{HV2}). 

As a consequence of definition~\ref{M-leg-submfld-conicpts}, $\Gh$ is a
compact manifold with corners in 
\begin{equation}
[\sfT_{\mf} X; \ \{ x = 0, \mu_1 = \mu_2 = \tau_1 = 0 \}],
\ilabel{Ghat-2}\end{equation}
with one boundary hypersurface
at the intersection of $\Gh$ and $\sfT _{\mf \cap H_i} M$ for each $i$ 
(which are
mutually nonintersecting, since the $H_i$ are mutually nonintersecting), and
one at the intersection of $\Gh$ and the front face of the blowup in
\eqref{Ghat-2}. If $q$ lies in the interior of $\Gh$  then the situation is as
for Legendrians in the scattering setup. If $q$ is on the boundary of $\Gh$,
but does not lie over $H_i$ for some $i$, then the situation is as for
Legendrian conic pairs as in the previous subsection. If $q$ is on the
boundary of $\Gh$ but not in $\Gs$ then the situation is as above. Thus the
only situation left to describe is if $q$ is in the corner of $\Gh$, lying
above $\ybar_0 \in \mf \cap H_i$ say. 

A local parametrization of $(G, \Gs)$ near $q$ (in coordinates as chosen
above) is a function 
$$\phi(x_1, y_1, y_2, s, w) = 1 + sx_1 \psi(x_1, y_1, y_2, s, w),
$$
with $\psi$ defined in a neighbourhood of $q' = (\ybar_0, 0, w_0) \subset M
\times [0, \infty) \times \RR^k$, such that 
$d_w \psi = 0$ at $q'$, 
$$
d_{y_2} \psi \ \text{ and } d_{(y_2,w)} \left( \frac{\pa \psi}{\pa
  w_i} \right) \text{ are linearly independent at } q',
$$
and such that near $q \in \Gh$, 
$$
\Gh = \{ (\ybar, -\phi, \frac{d_{x_1}(x_1 \psi)}{|d_{y_2} \psi|},
\frac{d_{y_1}(x_1 \psi)}{|d_{y_2} \psi|}, \widehat{ d_{y_2} \psi },
{|d_{y_2} \psi|} \},
$$
in the coordinates \eqref{Ghat-coords}. 

A Legendre distribution of order $(m,p; r_1, \dots r_d)$ associated to $(G,
\Gs)$ is a half-density 
such that for any $\upsilon_i \in \CI(M)$ whose support does not intersect
$H_k$, for $k \neq i$, $\upsilon_i u$ is of
the form $u = u_0 + (\sum_{j=1}^N u_j + \sum_{j=1}^M u_j')\nu$, where 
$u_0 \in \CIdot(X; \sfOh)$, and $u_j$, $u'_j$ have expressions 
\begin{equation}\begin{gathered}\label{eq:dist-1p}
u_j(x_1, x_2,y_1, y_2)
= \int_0^\infty ds \int dw \, e^{i\phi_j(x_1, y_1, y_2,v,w)/x} 
a_j(x_1, s, x_2/s, y_1, y_2,w)\\
\left( \frac{x_2}{s}\right)^{m-(k'+1)/2+N/4} s^{p-1+N/4}
x_1^{r_i-k/2+N/4-f_i/2} 
\end{gathered}\end{equation}
with $N=\dim M$, $a_j\in\Cinf_c([0,\epsilon)\times
U\times\RR^{k+k'})$, $U$ open in $\mf$, $f_i$ the 
dimension of the fibres of $H_i$ and
$\phi_j$ a phase function parametrizing a Legendrian $G$ on
$U$, and where $u'_j$ is as in \eqref{eq:dist-2}.

\subsection{The b-double space} Here we analyze the b-double space $\XXb$,
where $X$ is a compact manifold with boundary, from the perspective of
manifolds with corners with fibred boundaries. The manifold with corners
$\XXb$ has three boundary hypersurfaces: lb and rb, which are the lifts of
the left and right boundaries $\pa X \times X$, $X \times \pa X$ of $X^2$ to
$\XXb$, and bf, coming from the blowup of $(\pa X)^2$ (see figure
\ref{fig:XXb}). Thus, lb 
and rb have natural projections to $\pa X$. The fibres of lb and rb meet bf
transversally, so we may identify bf 
as the `main face' mf of $\XXb$. Given coordinates $(x,y)$ or
$z$ on $X$, we denote the lift to $\XXb$ via the left, resp. right
projection by $(x',y')$ or $z'$, resp. $(x'',y'')$ or $z''$. 
We may take the distinguished total boundary defining function to
be $x'$, for $\sigma = x'/x'' < C$ and $x''$, for $\sigma >
C^{-1}$. These are compatible since their ratio is constant on fibres on the
overlap region $C^{-1} < \sigma < C$ (this is trivially true since the
fibres of bf are points). 

These data
give $\XXb$ the structure of a manifold with corners with fibred boundary as
defined above. The $\sf$-vector fields then are the same as the sum of the scattering
Lie Algebra $\Vb(X)$ lifted to $\XXb$ from the left and right factors. 

On $\XXb$, and $\sfT \XXb$, it is most convenient to use coordinates lifted
from $X$ and $\sct X$. Near lb, 
but away from bf, we use coordinates $(x',y',z'')$ and coordinates $(\tau',
\mu', \zeta'')$ on $\sfT \XXb$ where we write a covector $q \in \sfT \XXb$
as
$$
q = \tau' \frac{dx'}{{x'}^2} + \mu' \cdot \frac{dy'}{x'} + \zeta'' \cdot
dz''.
$$
Similarly near rb, but away from bf, we use coordinates $(z',x'',y'';
\zeta',\tau'',\mu'')$. Near $\lb \cap \operatorname{bf}$, we use $(x'',
\sigma, y', y'')$ 
with corresponding coordinates $(\tau, \kappa, \mu', \mu'')$, by
writing $q \in \sfT \XXb$ as
$$
q = \tau \frac{dx'}{{x'}^2} + \kappa \frac{ d\sigma}{x'} + \mu'
\cdot \frac{dy'}{x'} + \mu'' \cdot \frac{dy''}{x''}.
$$
However, we may also use scattering cotangent coordinates $(\tau', \tau'', \mu',
\mu'')$ lifted from $\sct X$, where we write
$$
q = \tau' \frac{dx'}{{x'}^2} + \tau'' \frac{dx''}{{x''}^2} + \mu'
\cdot \frac{dy'}{x'} + \mu'' \cdot \frac{dy''}{x''}.
$$
This gives
\begin{equation}
\tau' = \tau + \sigma \kappa \quad \tau''
  = -\kappa.
\ilabel{tau}\end{equation}
The coordinates $(x'', \sigma, y', y'')$ hold good near bf as long as we
stay away from rb, when we need to switch to $(x', \sigma^{-1}, y', y'')$. 
The cotangent coordinates $(\tau', \tau'', \mu', \mu'')$ are good coordinates
globally near bf; notice that the roles of $(\tau, \tau_1, \mu_1, \mu_2)$
are played 
by $(\tau', \tau'', \mu', \mu'')$ near lb and $(\tau'', \tau', \mu'', \mu')$
near rb. 

The operator $H$ can act on half-densities on $\XXb$ by acting either on the
left or the right factor of $X$; these operators are denoted $H_l$ and $H_r$
respectively. For $H = \Lap + P$, where $P \in x \Diffsc^1(X)$, the Hamilton
vector field induced by $H_l$ and
the contact structure on $\sfT_{\bff} \XXbt$, with
respect to $x'$, takes the form
\begin{equation}
V_l = 2 \tau' \sigma \frac{\pa}{\pa \sigma} + 
2\tau' \mu' \frac{\pa}{\pa \mu'} - h' \frac{\pa}{\pa \tau'} + \left
( \frac{\pa h'}{\pa \mu'} \frac{\pa}{\pa y'} - \frac{\pa
h'}{\pa y'} \frac{\pa}{\pa \mu'}  \right) \quad h' = h(y', \mu').
\ilabel{Vl}\end{equation}
Similarly, the Hamilton vector field induced by $H_r$ and
the contact structure on $\sfT_{\bff} \XXbt$, with
respect to $x''$, takes the form
\begin{equation}
V_r = -2 \tau'' \sigma \frac{\pa}{\pa \sigma} + 
2\tau'' \mu'' \frac{\pa}{\pa \mu''} - h'' \frac{\pa}{\pa
  \tau''} + \left 
( \frac{\pa h''}{\pa \mu''} \frac{\pa}{\pa y''} - \frac{\pa h''}{\pa y''} 
\frac{\pa}{\pa \mu''}  \right) \quad h'' = h(y'', \mu'').
\ilabel{Vr}\end{equation}
Notice that $V_l$ and $V_r$ commute.

\section{Symbol calculus for Legendre distributions}
\subsection{Manifolds with boundary}
Let $X$ be a manifold with boundary of dimension $N$, and let
$$
u = x^q (2\pi)^{-k/2-n/4}
\left( \int e^{i\phi(y,v)/x} a(x,y,v) \, dv \right) \schd \in I^m(X, G; \scOh)
$$
be a Legendre distribution of order $m$. Let $C = \{ (y,v) \mid d_v \phi = 0
\}$ and let $\lambda$ be a set of functions in $(y,v)$-space such that
$(\lambda, d_v \phi)$ form local coordinates near $C$. We temporarily define 
the symbol relative to the
coordinate system $\mathcal{Z} = (x,y)$ and the parametrization $\phi$ to be
the half density on $G$ given by
\begin{equation}
\symrel{m}{Z}{\phi}(u) = \big( a(0,y,v) \restriction C \big) \left| \frac
{ \pa (d_v \phi, \lam)}{\pa (y,v)} \right| ^{-\half} |d\lam|^\half.
\ilabel{symrel}\end{equation}
Here we have used the correspondence \eqref{correspondence} between 
$C$ and $G$.

If we change coordinate system, the symbol changes by
\begin{equation}
\sigma^{m}_{\mathcal{Z}', \phi}(u) = \symrel{m}{Z}{\phi}(u)
a^{n/4-m}
e^{-i\rho(\mathcal{Z}', \mathcal{Z})}, \quad a = \frac{x'}{x}
\ilabel{change}\end{equation}
where
\begin{equation}
\rho(\mathcal{Z}', \mathcal{Z}) = \Big\{ a \mu_i \frac{\pa y_i}{\pa x'} 
 - \tau \frac{\pa a}{\pa x'} \Big\}
\restriction x = 0.
\ilabel{E-factor}\end{equation}
If the parametrization is changed, then by \cite{FIO1}, the symbol changes by
\begin{equation}
\symrel{m}{Z}{\psi}(u) = \symrel{m}{Z}{\phi}(u) e^{i\pi(\sign d^2_{vv}\psi -
  \sign d^2_{v'v'}\phi)/4};
\ilabel{M-factor}\end{equation}
the exponential is a locally constant function. We use these
transformation factors to define two line bundles, the $E$-bundle over $\sct
X$ which is
defined by the transition functions \eqref{E-factor}, and the Maslov bundle
over $G$ which is defined by the transition functions
\eqref{M-factor}. (These bundles will be described in much more detail in
\cite{HV4}.) Defining the bundle $\Sm (G) = |N^* \pa X|^{m-n/4} \otimes E
\otimes M(G)$ over $G$, we obtain an invariant symbol map from \eqref{symrel}
$$
\sigma^m : I^m(X, G; \scOh) \to \CI(G; \Omega^\half \otimes \Sm (G)).
$$

The elements of the symbol calculus for Legendre distributions on
manifolds with boundary have been given by Melrose and Zworski \cite{RBMZw}:

\begin{prop} The symbol map induces an exact sequence
$$
0 \to I^{m+1}(X, G; \scOh) \to I^m(X, G; \scOh) \to
\CI(G, \Omega^\half \otimes S^{[m]}(G)) \to 0.
$$
If $P \in \Psi^k(X; \scOh)$ and $u \in I^m(X, G; \scOh)$, 
then $Pu \in I^m(X, G; \scOh)$ and 
$$
\sigma^m(Pu) = \Big( \sigma(P) \rest G \Big) \sigma^m(u).
$$
Thus, if the symbol of $P$ vanishes on $G$, then $Pu \in I^{m+1}(X, G;
\scOh)$. The symbol of order $m+1$ of $Pu$ in this case is
\begin{equation}
\Big( -i \mathcal{L}_{H_p} -i \big(\half + m - \frac{N}{4} \big) \frac{ \pa
  p}{\pa \tau} + p_{\sub} \Big) \sigma^m(u) \otimes |dx|,
\ilabel{transport}\end{equation}
where $H_p$ is the Hamilton vector field of $p$, the principal symbol of
$P$, and $p_{\sub}$ is the subprincipal symbol of $P$. 
\end{prop}

The symbol calculus for intersecting Legendre distributions is easily
deduced from Melrose and Uhlmann's calculus of intersecting Lagrangian
distributions. The symbol takes values in a bundle over $L_0 \cup L_1$. Let
$\rho_1$ be a boundary defining function for $\pa L_1$ as a submanifold of
$L_0$, and $\rho_0$ be a boundary defining function for $\pa L_1$ as a
submanifold of $L_1$. To
define the symbol, note that the symbol on $L_0$ is defined by continuity from
distributions microsupported away from $L_1$, and takes values in
\begin{equation}
\rho_1^{-1} \CI(\Omega^{1/2}(L_0)\otimes S^{[m+1/2]}(L_0)) = 
\rho_1^{-1/2} \CI(\Omega_b^{1/2}(L_0 \setminus \pa L_1)
\otimes S^{[m+1/2]}(L_0)),
\ilabel{int-Leg-dist}\end{equation}
while the symbol on
$L_1$ defined by continuity from
distributions microsupported away from $\pa L_1$ takes values in
$$
\CI(\Omega^{1/2}(L_1) \otimes S^{[m]}(L_1)) =
\rho_0^{1/2}\CI(\Omega_b^{1/2}(L_1) \otimes S^{[m]}(L_1)).
$$
Melrose and Uhlmann showed
that the Maslov factors were canonically isomorphic on $L_0 \cap L_1$, so
$S^{[m]}(L_0)$ is naturally isomorphic to $S^{[m]}(L_1)$ over $L_0 \cap
L_1$. Canonical restriction of the
half-density factors to $L_0 \cap L_1$ gives terms in $\CI(\Omega^\half(L_0
\cap L_1) \otimes S^{[m]}(L_1) \otimes |N^*_{L_0}\pa L_1|^{-1/2} \otimes
|N^* \pa X|^{1/2}$ and $\CI(\Omega^\half(L_0
\cap L_1) \otimes S^{[m]}(L_1) \otimes |N^*_{L_1}\pa L_1|^{1/2}$
respectively. In 
fact $|N^*_{L_0}\pa L_1| \otimes |N^*_{L_1}\pa L_1| \otimes |N^* \pa
X|^{-1}$ is canonically trivial; an
explicit trivialization
is given by 
\begin{equation}
(d\rho_0, d\rho_1, x^{-1}) \mapsto x^{-1} \omega(V_{\rho_0},
V_{\rho_1}) \restriction L_0 \cap L_1,
\ilabel{triv}\end{equation}
where $V_{\rho_i}$ are the Hamilton vector fields of the
functions $\rho_i$ extended into $\sct X$, and $\omega$ is the standard
symplectic form. Thus the two bundles are naturally isomorphic over the
intersection. We define the bundle $S^{[m]}(\Lt)$ to be that bundle such
that smooth sections of $\Omega^{1/2}_b(\Lt) \otimes S^{[m]}(\Lt)$ are precisely
those 
pairs $(a,b)$ of sections of $\rho_1^{-1} \CI(\Omega^{1/2}(L_0)\otimes
S^{[m+1/2]}(L_0))$ and $\rho_0^{1/2}\CI(\Omega^{1/2}_b(L_1) \otimes 
S^{[m]}(L_1))$ such that 
\begin{equation}
\rho_1^{1/2} b = e^{i\pi/4} (2\pi)^{1/4} \rho_0^{-1/2} a \text{ at }L_0 \cap L_1
\ilabel{compat}\end{equation}
under the
identification \eqref{triv}. The symbol maps of order $m$ on $L_1$ and
$m+1/2$ on $L_0$ then extend in a natural way to a symbol map of order $m$
on $\Lt$ taking values in $\Omega_b^{1/2}(\Lt) \otimes S^{[m]}(\Lt)$. 

\begin{prop} The symbol map on $\Lt$ yields an exact sequence
\begin{equation}
0 \to I^{m+1}(X, \Lt; \scOh) \to I^m(X, \Lt; \scOh) \to
\CI(\Lt, \Omega^\half_b \otimes S^{[m]}) \to 0. 
\ilabel{ex-int-1}\end{equation}
Moreover, if we consider just the symbol map to $L_1$, there is an exact
sequence
\begin{multline}
0 \to I^{m+1}(X, \Lt; \scOh) + I^{m+\half}(X, L_0; \scOh) \to I^m(X, \Lt;
\scOh)  \\  
\to \CI(L_1, \Omega^\half \otimes S^{[m]}) \to 0.
\ilabel{ex-int-2}\end{multline}
If $P \in \Psi^k(X; \scOh)$ and $u \in I^m(X, \Lt; \scOh)$, 
then $Pu \in I^m(X, \Lt; \scOh)$ and 
$$
\sigma^m(Pu) = \Big( \sigma(P) \rest \Lt \Big) \sigma^m(u).
$$
Thus, if the symbol of $P$ vanishes on $L_1$, then $Pu$ is an element of
$I^{m+1}(X, \Lt;
\scOh) + I^{m}(X, L_0; \scOh)$. The symbol of order $m+1$ of $Pu$ on $L_1$
in this case is given by \eqref{transport}. 
\end{prop}

For a conic pair of Legendre submanifolds $\Gt = (G, \Gs)$, with $\Gh$ the
desingularized submanifold obtained by blowing up $\Gs$, the symbol is
defined by continuity from the regular part of $G$. The symbol calculus then
takes the form

\begin{prop} Let $s$ be a boundary defining function for $\Gh$. Then 
there is an exact sequence
$$
0 \to I^{m+1,p}(X, \Gt; \scOh) \to I^m(X, \Gt; \scOh) \to
s^{m-p}\CI(\Gh, \Omega^\half_b \otimes S^{[m]}(\Gh)) \to 0.
$$
If $P \in \Psi^k(X; \scOh)$, and $u \in I^{m,p}(X, \Gt; \scOh)$, 
then $Pu \in I^{m,p}(X, \Gt; \scOh)$ and 
$$
\sigma^m(Pu) = \Big( \sigma(P) \rest \Gh \Big) \sigma^m(u).
$$
If the symbol of $P$ vanishes on $G$, then $Pu \in I^{m+1,p}(X, \Gt;
\scOh)$. The symbol of order $m+1$ of $Pu$ in this case is given by
\eqref{transport}. 
\end{prop}

\subsection{Codimension two corners} When we have codimension two corners,
then essentially the same results hold by continuity from the main face. 
The symbol is defined as a half-density on $G$ by continuity from the
interior of mf, where the scattering situation applies. 
We must restrict to
differential operators, however, since pseudodifferential operators have not
been defined in this context. 

Let $M$ be a manifold with codimension 2 corners with fibred boundaries,
let $N = \dim M$, and
let $G$ be a Legendre distribution. Let $\rho_i$ be a boundary defining 
function for $H_i$. The Maslov bundle $M$ and the E-bundle are defined via
the scattering structure over the interior of $G$ and extend to smooth
bundles over the whole of $G$ (that is, they are smooth up to each boundary
of $G$ at $\sfT _{H_i \cap \mf} M$). Let $S^{[m]}(G) = M(G) \otimes E
\otimes |N^* \mf|^{m-N/4} \otimes |N^* H_1|^{m-N/4} \otimes \dots \otimes
|N^* H_d|^{m-N/4}$. Finally let ${\bf r}$ stand for $(r_1, \dots, r_d)$,
and let $\rho^{\bf r} = \prod_i \rho_i^{r_i}$. 

\begin{prop} There is an exact sequence
$$
0 \to I^{m+1,{\bf r}}(M, G; \sfOh) \to I^{m,{\bf r}}(M, G; \sfOh) \to 
\rho^{m-{\bf r}} \CI(G, \Omega^\half_b \otimes S^{[m]}(G)) \to 0.
$$
If $P \in \Diff(M; \sfOh)$ and $u \in I^{m,{\bf r}}(M, G; \sfOh)$, 
then $Pu \in I^{m,{\bf r}}(M, G; \sfOh)$ and 
$$
\sigma^m(Pu) = \Big( \sigma(P) \rest G \Big) \sigma^m(u).
$$
Thus, if the symbol of $P$ vanishes on $G$, then $Pu \in I^{m+1,{\bf r}}(M, G;
\sfOh)$. The symbol of order $m+1$ of $Pu$ in this case is given by
\eqref{transport}. 
\end{prop}

For a conic pair of Legendre submanifolds $\Gt = (G, \Gs)$, with $\Gh$ the
desingularized submanifold obtained by blowing up $\Gs$, the symbol calculus
takes the form

\begin{prop}\label{ex-conic-2} Let $s$ be a boundary defining function for
$\Gh$ at $\Gh \cap \Gs$. Then  
there is an exact sequence
\begin{multline}
0 \to I^{m+1,p; {\bf r}}(M, \Gt; \sfOh) \to I^{m,p; {\bf r}}(M, \Gt; \sfOh) \\
\to \rho^{m-{\bf r}} s^{m-p} \CI(\Gh, \Omega^\half_b \otimes S^{[m]}(\Gh))
\to 0. 
\end{multline}
If $P \in \Diff(M; \sfOh)$, and $u \in I^{m,p;{\bf r}}(M, \Gt; \sfOh)$, 
then $Pu \in I^{m,p;{\bf r}}(M, \Gt; \sfOh)$ and 
$$
\sigma^m(Pu) = \Big( \sigma(P) \rest \Gh \Big) \sigma^m(u).
$$
If the symbol of $P$ vanishes on $G$, then $Pu \in I^{m+1,p;{\bf r}}(M, \Gt;
\sfOh)$. The symbol of order $m+1$ of $Pu$ in this case is given by
\eqref{transport}. 
\end{prop}

The proofs of these propositions are omitted, since they are easily deduced
from the codimension one case.

\section{Parametrix construction}
In this section, we consider self-adjoint operators $H$ of the form $\Lap +
P$, where $\Lap$ is the positive Laplacian with respect to a short-range
metric on a compact manifold with boundary, $X$, and $P \in x^2
\Diffsc^1(X)$ is a short-range perturbation of $\Lap$. In the following
section, we consider metrics and perturbations of long range gravitational
type. Let $R(\sigma)$ denote the resolvent $(H - \sigma)^{-1}$ of $H$. 

In this section, we 
directly construct a parametrix $G(\lam)$ for \eqref{eq-Rt} whose error term
$E(\lam) = (H - \lam^2) 
G(\lam) - \Id$ is compact. Using Fredholm theory and a unique continuation
theorem we solve away the error, giving us a Schwartz kernel $\Rt(\lam)$.
We then show that $\Rt(\sqrt{\sigma})$ has an analytic continuation (as a
distribution on $X^2$) 
to the upper half $\sigma$ plane which
agrees with the resolvent $R(\sigma)$ there. This proves that $\Rt(\lam)$ and
$R(\lam^2 + i0)$ coincide on the real axis. 

The distribution $\Rt(\lam)$ has the defining property that
\begin{equation}
(H - \lam^2) \Rt(\lam) = \Id,
\ilabel{defining}\end{equation}
as an operator on $\CIdot(X; \scOh)$, and that
\begin{equation}
\scwf_{\XXbt}(\Rt(\lam)) \subset K_- \text{ as defined in } \eqref{Kminus}.
\ilabel{outgoing}
\end{equation}
Here $\XXbt = \XXb \setminus \{ \lb \cup \rb \}$ is regarded as an open
manifold with boundary, so that we can talk about the scattering wavefront
set over the interior of bf. Equation \eqref{outgoing} is the microlocal
version of the outgoing Sommerfeld radiation condition. (For example, if
$\lam > 0$, 
$e^{i\lam/x}$ has wavefront set in $K_-$, while $e^{-i\lam/x}$ does not.)

Equation \eqref{defining} means that the kernel of $\Rt(\lam)$, which we also
denote by $\Rt(\lam)$ by an abuse of notation, satisfies
\begin{equation}
(H_l - \lam^2) \Rt(\lam) =K_{\Id}, 
\ilabel{diag}\end{equation}
where $K_{\Id}$ is the kernel of the identity operator, i.e.\ it is
a delta distribution on the diagonal, and $H_l$ is the operator
$H$ acting on the left factor of $X$ in $X\times X$.

There are four main steps in the construction. First we find an
approximation to $\Rt(\lam)$ in the scattering calculus, $G_1(\lam)
\in \Psisc^{-2}(X)$, which removes the singularity on the diagonal in
\eqref{diag}. This leaves an error which, when viewed on the b-double space
$\XXb$, is singular at the boundary of the diagonal $\pa \diag_\bl$. In
fact, it is Legendrian at a Legendre submanifold lying over $\pa \diag_\bl$
which we denote $N^* \diag_\bl$ (see \eqref{Ndiag}). We solve this error
away locally 
near $\pa \diag_\bl$ using an intersecting Legendrian construction which is
due (in the Lagrangian setting) to Melrose and Uhlmann
\cite{Rbm-Uhl:Intersecting}; the 
singularities inside $N^* \diag_\bl \cap \Sigma(H_l - \lam^2)$ propagate in a
Legendre submanifold $\Lpl$. This Legendre submanifold intersects both lb
and rb, and an `outgoing' Legendre submanifold $\Lsl$; $(\Lpl, \Lsl)$ form
an conic pair of Legendre submanifolds and we can find a conic Legendre pair
which solves away the error up to an error term which is Legendrian only at
$\Lsl$, ie we can solve away the errors at $\Lpl$ completely. Finally, this
outgoing error is solved away, using a very standard argument in scattering
theory, at lb and bf, leaving an error $E(\lam)$ which is compact on
weighted $L^2$ spaces $x^l L^2(X)$ for all $l > 1/2$. 

Thus, we seek $G(\lam)$ (and $\Rt(\lam)$) in the class 
\begin{equation}
\Psisc^{2,0}(X) + I^{-\frac1{2}}(N^* \diag_\bl, \Lpl; \scOh) +
I^{-\frac{1}{2}, 
  \frac{n-2}{2}; \frac{n-1}{2}, \frac{n-1}{2}}(\Lpl, \Ls(\lam); \scOh)
\ilabel{seek}\end{equation}
where the second term is an intersecting Legendrian distribution and the
third is a Legendre conic pair with orders $-1/2$ at $\Lpl$, $(n-2)/2$ at
$\Lsl$ and $(n-1)/2$ at lb and rb. In this class of distributions there is a
unique solution $\Rt(\lam)$ to \eqref{defining} and \eqref{outgoing}. 

To avoid cumbersome notation, $Q$ will denote a generic correction to
the para\discretionary{-}{metrix}{metrix}
constructed so far, and $E$ will denote a generic error. The values of these
symbols is allowed to change from line to line.

\subsection{Pseudodifferential approximation}\label{sec:pseudodiff-approx}
The first step in constructing
$G(\lam)$ is a very 
standard argument. We seek $G_1(\lam) \in \Psisc^{-2}(X)$ such that
$$
(H - \lam^2) G_1(\lam) = \Id + E_1(\lam), \quad E_1(\lam) \in
\Psisc^{-\infty}(X).
$$
This will mean that the error term $E_1(\lam)$ has a smooth kernel
(times the standard half-density) on $\XXsc$, so that we have solved away
completely 
the singularity along the diagonal. 

The standard elliptic argument applies here 
since the interior symbol of $H - \lam^2$
is elliptic. Thus, we first choose any $Q \in \Psisc^{-2}(X)$ whose
interior symbol is $|\cdot|_g^{-2} = (\sigma^2(H - \lam^2))^{-1}$. Then
$$
(H - \lam^2) Q = \Id + E, \quad E \in \Psisc^{-1}(X).
$$
Multiplying $Q$ by a finite Neumann series $(\Id + E + \dots + E^{k-1})$ thus
gives an error $E^k \in \Psisc^{-k}(X)$. Taking an 
asymptotic limit gives us a $G_1(\lam) \in \Psisc^{-2}(X)$ with the
desired error term.

\subsection{Intersecting Legendrian construction} 
\label{sec:intersecting-Leg-constr}

In the next step of the construction we move to $\XXb$, and view the error
$E_1(\lam)$ from the first step of the construction on $\XXb$ rather than
$\XXsc$. On $\XXb$ it has a smooth kernel except at $\pa \diag_\bl$ where it
has a conic singularity. That is, at $\pa \diag_\bl = \{ x' = 0,
\sigma \equiv x'/x'' = 1, y' = y'' \}$,  the kernel is a smooth 
(and compactly supported) function of $x'$, $S = (\sigma-1)/x'$, $Y =
(y'-y'')/x'$ and $y'$; this is easy to see since these are smooth
coordinates on sf $\subset \XXsc$. Using the Fourier transform, we write
$$
E_1(\lam) = \left( \int_{\RR^n} e^{i((y'-y'')\cdot \eta + (\sigma - 1)
  t)/x'} 
a(x',y',\eta, t) \, d\eta \, dt \right) \nu 
$$
where $a$ is smooth in all variables, and in addition Schwartz in $(\eta,
t)$. The phase function $(y'-y'')\cdot \eta + (\sigma - 1)t$ parametrizes
the Legendrian
\begin{equation}
N^* \diag_\bl = \{ y' = y'', \sigma = 1, \mu' = -\mu'', \tau' = -\tau'' \}.
\ilabel{Ndiag}\end{equation}
Therefore, $E_1(\lam)$ is a Legendre distribution of order $0$
associated to $N^* \diag_\bl$. (To be pedantic, $E_1(\lam)$
does not fall strictly in the class of Legendre distributions as defined by
Melrose and Zworski since its microsupport is not compact; from
\eqref{Ndiag} we see that the microsupport is a vector
bundle over $\pa \diag_\bl$. It is instead an `extended Legendre distribution'
as defined in \cite{Hassell:Plane}. However this is of no significance since
the symbol 
is rapidly decreasing in each fibre of the vector bundle, hence all
constructions we wish to perform here are valid in this context.)

Observe that $\symbb(H_l - \lam^2) = {\tau'}^2 + h(y',\mu') - \lam^2$ vanishes
on a codimension one submanifold of $N^* \diag_\bl$, and does so
simply. Consider the vector field $V_l$ which is given by \eqref{Vl}. Since
${\tau'}^2 
+ h = \lam^2 \neq 0$ on $\Sigma(H_l - \lam^2)$, at least one of the coefficients
of $\pa_\sigma$ and $\pa_{\tau'}$ in \eqref{Vl} is nonzero, so $V_l$ is
transverse to $N^* 
\diag_\bl$ at the intersection with $\Sigma(H_l - \lam^2)$. We define
$L^\circ(\lam)$ to be the flowout Legendrian from $N^* \diag_\bl \cap
\Sigma(H_l - \lam^2)$ 
with respect to $V_l$, and $L^\circ_\pm(\lam)$ to be the flowout in the
positive, resp. negative direction with respect to $V_l$. Thus, at least
locally near $N^* \diag_\bl$, $L^\circ_\pm (\lam)$ are smooth manifolds with
boundary. Notice that by \eqref{tau}, $N^* \diag_\bl$ is contained in $\tau
= 0$ and $V_l(\tau) < 0$ at $N^* \diag_\bl$. Thus, at least locally near
$N^* \diag_\bl$, $L^\circ_+(\lam)$ is contained in $K_- = \{ \tau \leq 0 \}$.
The global properties of $L^\circ_\pm (\lam)$ are studied in the next
section; in this section we only work microlocally near $N^* \diag_\bl$.

The first step in solving away the error $E_1(\lam)$ from the previous step is
to find an intersecting Legendrian $Q \in 
I^{-1/2} (\XXb, (N^*\diag_\bl, \Lpl); \sfOh)$ such that 
\begin{equation}
(H_l - \lam^2) Q - E_1(\lam) \in I^{1}(\XXb, N^* \diag_\bl; \sfOh) +
I^{\frac{3}{2}} (\XXb, (N^* \diag_\bl, \Lpl), \sfOh),
\ilabel{err-1}\end{equation}
microlocally near $N^* \diag_\bl$. To do this we choose $Q$ with symbol on
$N^* \diag_\bl$ equal to $\symbb^0(H_l - \lam^2)^{-1}
\sigma^0(E_+^{(1)}(\lam))$. This is an admissible symbol on $N^* \diag_\bl$ by
\eqref{int-Leg-dist}, and \eqref{ex-int-1}, 
since $\symbb^0(H_l - \lam^2)$ is a boundary defining
function for $\Lpl$ on $N^* \diag_\bl$. It determines the value of the
symbol on $\Lpl$ at the boundary by \eqref{compat}. We extend this symbol by
requiring that the transport equation, \eqref{transport}, be satisfied. 
This equation is a first order
linear ODE with smooth coefficients, so there is a unique solution in a
neighbourhood of $N^* \diag_\bl$. Then the symbol of order $-1/2$ of 
$(H_l - \lam^2) R_0 - E_+^{(1)}(\lam)$ vanishes, and there is an additional
order of vanishing on $\Lpl$ since the transport equation is satisfied.
Thus by \eqref{ex-int-2} the error term is as in \eqref{err-1}.

We now show inductively that we can solve away an error $E_k$ which is in
$I^{k}(N^* \diag_\bl) + I^{k+1/2} (N^*\diag_\bl, \Lpl)$ with a term
$Q_k \in I^{k+1/2}(N^*\diag_\bl, \Lpl)$, 
up to an error which is in $I^{k+1}(N^* \diag_\bl) + I^{k+3/2}
(N^*\diag_\bl, \Lpl)$. The argument is the same as above: we take the symbol
of order $k$ on $N^*\diag_\bl$ equal to $\symbb^0(H_l - \lam^2)^{-1}
\sigma^k(E_k)$, and the symbol on $\Lpl$ to solve away the symbol of
order $k+1/2$ of $E_k$ when the transport operator is applied to it. 
Taking an asymptotic sum of $Q$ and the $Q_k$'s gives us an error term
which is microlocally trivial near $N^* \diag_\bl$. 
By cutting off away from $\pa \diag_\bl$, we obtain an error 
$$
E_2(\lam) \in I^{1/2}_c(\Lpl),
$$
where the subscript $c$ indicates that the microlocal support is compact and
disjoint from the intersection with $N^* \diag_\bl$.

\subsection{Structure of $L(\lam)$}
In this section we analyze the global structure of $L^\circ(\lam)$. This was
defined as the flowout from $N^* \diag_\bl \cap \Sigma(H_l - \lam^2)$ by the
vector field $V_l$. In fact, it is quite easy to see that $N^* \diag_\bl
\cap \Sigma(H_l - \lam^2) = N^* \diag_\bl \cap \Sigma(H_r - \lam^2)$. Moreover,
neither $V_l$ nor $V_r$ is tangent to $N^* \diag_\bl$ at any point contained
in $\Sigma(H_l - \lam^2)$, but the difference $V_l - V_r$ is tangent to $N^*
\diag_\bl$. Since $V_l$ and $V_r$ commute, this shows that the flowout with
respect to $V_l$ is the same as the flowout with respect to $V_r$. We will
soon see that the symbols of our parametrix on $\Lpl$, defined so as to
satisfy the left transport equation, also satisfy the right transport
equation.

It is convenient to write down $L^\circ(\lam)$ explicitly. Indeed, the
computation of Melrose and Zworski can be applied with a minor change
(that takes care of the behavior in $\sigma$) to deduce that
\begin{equation}\begin{split}\label{eq:sp-1c}
L^\circ(\lam) 
=&\{(\theta,y',y'',\tau',\tau'',\mu',\mu''):\ \exists(y,\muh)\in S^*\bX,
\ s,s'\in(0,\pi),\ \text{s.t.}\\
&\quad \sigma = \tan \theta = \frac{\sin s'}{\sin s},
\ \tau'=\lambda\cos s',
\ \tau''=-\lambda\cos s,\\
\quad(y',\mu')&=\lambda\sin s'\exp(s'H_{\half h})(y,\muh),
(y'',\mu'')=-\lambda\sin s \exp(s H_{\half h})(y,\muh)\}\\
&\cup T_+(\lam) \cup T_-(\lam), \quad T_\pm (\lam) = 
\{(\sigma,y,y,\pm \lambda,\mp \lambda,0,0):\ \sigma\in(0,\infty),\ y\in\bX\}
\end{split}\end{equation}
The sets $T_\pm (\lam)$ are, 
for fixed $y$, integral curves of both vector fields, and they appear
separately only because we used the parameterization of Melrose-Zworski.
The smooth structure near $T_\pm (\lam)$ follows from the
flowout description, but is not apparent in this parameterization;
we discuss it below while describing the closure of $L^\circ(\lam)$.

The closure $L(\lam)$ of $L^\circ(\lam)$ is $\sfT\XXb$ is clear from the
above description; it is
\begin{equation}
\cl L =L(\lam)\cup F_{\lambda}\cup F_{-\lambda}
\end{equation}
where
\begin{equation}\begin{split}\label{eq:sp-1}
L(\lam) 
=&\{(\theta,y',y'',\tau',\tau'',\mu',\mu''):\ \exists(y,\muh)\in S^*\bX,
\ s,s'\in[0,\pi],\\
&\quad (\sin s)^2+(\sin s')^2>0,\ \text{s.t.}\\
&\quad \sigma = \tan \theta = \frac{\sin s'}{\sin s},
\ \tau'=\lambda\cos s',
\ \tau''=-\lambda\cos s,\\
\quad(y',\mu')&=\lambda\sin s'\exp(s'H_{\half h})(y,\muh),
(y'',\mu'')=-\lambda\sin s \exp(s H_{\half h})(y,\muh)\}\\
& \cup T_+(\lam) \cup T_-(\lam), \quad \text{and} \\
F_{\lam} &= \{(\sigma,y',y'',-\lam,-\lam,0,0) \mid 
\exists \text{ geodesic of length } \pi \text{ connecting } y', \, y'' \}.
\end{split}\end{equation}

Note that the requirement $(\sin s)^2+(\sin s')^2>0$ just
means that $s$ and $s'$ cannot take values in $\{0,\pi\}$ at the same time.
The set $L(\lam) \setminus L^o(\lam)$ comprises those points where one of $s$,
$s'$ takes values in $\{ 0, \pi \}$ while the other lies in $(0, \pi)$. The
sets $T_\pm (\lam)$ in \eqref{eq:sp-1} comprise the limit points where $s$ and
$s'$ converge either both to $0$ or both to $\pi$, whilst $F_{\pm \lam}$
comprise the limit points as $s \to 0$ and $s' \to \pi$ or vice versa.

The smooth structure near $T_\pm (\lam)$ becomes apparent 
if we note that near $\tau'=\lambda$,
$\tau''=-\lambda$, $\sigma\in[0,C)$ where $C>1$, $L(\lam)$ is
given by
\begin{equation}\begin{split}\label{eq:sp-1-a}
\{(\sigma&,y',y'',\tau',\tau'',\mu',\mu''):\ \exists(y,\mu)\in T^*\bX,
\ |\mu|<C^{-1},\ \sigma\in[0,C)\ \text{s.t.}\\
&\tau'=\lambda(1-|\sigma\mu|^2)^{1/2},
\ \tau''=-\lambda(1-|\mu|^2)^{1/2},\\
&(y',\mu')=\lambda\exp(f(\sigma\mu)V_h)(y,\sigma\mu),
\ (y'',\mu'')=-\lambda\exp(f(\mu)V_h)(y,\mu)\}
\end{split}\end{equation}
where
$f(\mu)={|\mu|}^{-1}\arcsin(|\mu|)$
is smooth and $f(0)=1$. Thus, the differential of the map
\begin{equation}
(y,\mu)\mapsto-\lambda\exp(f(\mu)V_h)(y,\mu)=(y'',\mu'')
\end{equation}
is invertible near $\mu=0$, so it gives a diffeomorphism near $|\mu|=0$.
Hence, $\sigma$ and $(y'',\mu'')$ give coordinates on $L(\lam)$ in this
region, so $L(\lam)$ is smooth here. Away from $T_+(\lam)$, coordinates on
$L(\lam)$ are $\sigma$, $y''$, $\muh''$ and $s$. 

In the coodinates $(y, \hat \mu, s, s')$, the vector field $V_l$ is given by
$\sin s' \pa_{s'}$ and $V_r$ is given by $\sin s \pa_s$. The intersection of
$L(\lam)$ and $N^* \diag_\bl$ is given by $\{ s = s' \}$. Thus
$\Lpl$ is given by $\{ s \leq s' \}$. On $\Lpl$, $\tau = \tau' + \sigma
\tau''$ by \eqref{tau}, so 
$$
\tau = \lam \frac{\sin (s-s')}{\sin s} \leq 0 \text{ on } \Lpl.
$$
Thus, any distribution in $I^m(N^* \diag_\bl, \Lpl)$ satisfies condition
\eqref{outgoing}. 

We also define
\begin{equation}
\Ls(\lam)=\{(\theta,y',y'',-\lam,-\lam,0,0):\ y',y''\in\bX,\ \theta
\in[0,\pi/2]\},
\label{eq:L-sharp}\end{equation}
so $\Ls(\lam)$ is a Legendrian submanifold of $\sct_{\ffb}\XXbt$, and
\begin{equation}
\cl L\cap\Ls(\pm \lam)=F_{\pm\lambda}.
\end{equation}

\begin{prop}\ilabel{conicpair} The pair 
\begin{equation}
\Lt (\lam) = (L(\lam), \Ls(\lam) \cup \Ls(-\lam))
\end{equation}
is a pair of intersecting Legendre manifolds with conic points.
\end{prop}

\begin{proof} 
We must show that when the set $\{ t q \mid t > 0, q \in \Lsl \}$ is blown
up inside $\sfT \XXb$, the closure of $L(\lam)$ is a smooth manifold with
corners which meets the front face of the blowup transversally. 
Let us restrict attention to a neighbourhood of $\Lsl$; the case of
$\Ls(-\lam)$ is similar. 
Consider the vector field $V_l + V_r$. By \eqref{Vl} and \eqref{Vr}, in
$\Sigma(H_l - \lam^2) \cap \Sigma(H_r - \lam^2)$ this is given by
$$
2(\tau' - \tau'') \sigma \pa_\sigma + 2\tau' \mu' \cdot \pa_{\mu'} + 
2\tau'' \mu'' \cdot \pa_{\mu''} + \pa_{\mu'} h' \cdot \pa_{y'} - 
\pa_{y'} h' \cdot \pa_{\mu'} + \pa_{\mu''} h'' \cdot \pa_{y''} - 
\pa_{y''} h'' \cdot \pa_{\mu''}
$$
This is
equal to $-2\lam$ times the b-normal vector field $\mu' \cdot \pa_{\mu'} +
\mu'' \cdot \pa_{\mu''}$ plus a sum of vector fields which have the form 
$\rho V$, where $\rho$ vanishes at $\Lsl$ and $V$ is tangent to lb and
$\Lsl$ (all
considerations taking place inside $\Sigma(H_l - \lam^2) \cap \Sigma(H_r
- \lam^2)$). Thus, under blowup of $\{ t q \mid t > 0, q \in \Lsl \}$, $V_l +
V_r$ lifts to a vector field of the form
\begin{equation}
V_l + V_r = -2\lam s \pa_s + s W,
\ilabel{VlVr}\end{equation}
where $W$ is smooth and tangent to the boundary of $\hat L(\lam)$, 
and so dividing by $s$
yields a nonvanishing normal vector field plus a smooth tangent
vector field. As above, such a vector field has a continuation
across the boundary to the double of $\hat L(\lam)$ (across the front face) as
a smooth 
nonvanishing vector field. This holds true smoothly up to the corner with
lb, so $\hat L(\lam)$ is a smooth manifold with corners. 
\end{proof}

\subsection{Smoothness of symbols}\label{smoothness} In the next stage of
the construction, we 
solve away the error $E_2(\lam)$ which is microsupported in the interior
of $\Lpl$. This involves solving the transport equation globally on
$\Lpl$. In view of Proposition~\ref{conicpair}, we can expect the
construction to involve Legendrian conic pairs with respect to $(\Lpl,
\Lsl)$. In order for the symbol to be quantizable to such a conic pair, we
need to show regularity of the symbol on $\Lplh$, so that it lies in the
symbol space of the exact sequence from Proposition~\ref{ex-conic-2}. 
That is, the symbol of order $j-1/2$ on $\Lpl$ should lie in
\begin{equation}
\rho_{\lb}^{n/2-j} \rho_{\rb}^{n/2-j} \rho_{\#}^{(n-1)/2-j} \CI(\Lplh; 
\Omega_b^{1/2} \otimes S^{[j-1/2]}(\Lplh)).
\ilabel{regularity}
\end{equation}
(We will ignore the symbol bundle $S^{[j-1/2]}$ in the rest of this
section.)

To do this, we observe that the symbol on $\Lpl$ automatically satisfies
the transport equation for the right Hamilton vector field. To see this,
let $G_2(\lam)$ be the approximation to $\Rt(\lam)$ constructed so far, with
$$
(H_l - \lam^2) G_2(\lam) - K_{\Id} = E_2(\lam),
$$
and
$$
\scwf(E_2(\lam)) \subset \Lpl \cap \{ \tau < -c \}
$$
for some $c > 0$. Consider applying $H_r - \lam^2$ to $E_2(\lam)$. Since $H_l$
and $H_r$ commute, and $H_l K_{\Id} = H_r K_{\Id}$, we have
$$\begin{aligned}
(H_l - \lam^2) \left( (H_r - \lam^2) G_2(\lam) - \Id \right) 
&= (H_r - \lam^2) \left ( (H_l - \lam^2) G_2(\lam) - \Id \right) \\ 
&= (H_r - \lam^2) E_2(\lam).
\end{aligned}$$
We claim that $\scwf( (H_r - \lam^2) G_2(\lam) - K_{\Id})$ is contained in $
\{ \tau < -c \}$. For if there is a point where $\tau \geq -c$, then by
\eqref{prop-3}, the maximal bicharacteristic ray in 
\begin{multline}
\Sigma(H_l - \lam^2) \setminus \scwf\Big((H_l - \lam^2)\big((H_r - \lam^2) G_2(\lam)
- K_{\Id}\big)\Big) \\ = \Sigma(H_l - \lam^2) \setminus \scwf((H_r - \lam^2)
E_2(\lam))
\end{multline}
lies in $\scwf( (H_r - \lam^2) G_2(\lam) - K_{\Id})$. Such rays always propagate
into $\tau > 0$. But
$$
\scwf(G_2(\lam)) \subset \{ \tau \leq 0 \}, \quad
\scwf(K_{\Id}) \subset \{ \tau \leq 0 \}, 
$$
so this is impossible. Consquently, $(H_r - \lam^2) G_2(\lam) - K_{\Id}$ has
no scattering wavefront set
for $\{ \tau \geq -c \}$, and so the symbols of $G_2(\lam)$ must obey the
right transport equations in this region. By cutting off the symbols closer
and closer 
to the boundary of $\Lpl$, we see that the right transport equations must be
satisfied everywhere on $\Lpl$. 

Let us examine the form of these transport equations at the boundary of
$\Lpl$. Near lb, we have coordinates $(y'', \mu'', \sigma)$ near $T_+(\lam)$
and $(y'', \muh'', \sigma, s)$ away from $T_+(\lam)$, which are valid
coordinates for $\sigma < 2$, say. The situation near rb is similar so the
argument will be omitted. 

In either set of coordinates, the left vector field, restricted to $\Lpl$, takes
the form 
$$
V_l = 2\tau' \sigma \pa_\sigma .
$$
Also, by Lemma~\ref{lem:sub-pr}, the subprincipal symbol of $H_l - \lam^2$,
which is equal to the subprincipal symbol of $H - \lam^2$ in the
singly-primed coordinates, vanishes where $\mu'$ vanishes, and $\mu' = 0$ at
lb on $\Lpl$. Therefore, by \eqref{transport}, the transport equation for the
symbol of order $-1/2$ takes the
form
\begin{equation} 
\Big( -i\big({\mathcal L}_{V_l} 
- n \tau' \big) + 
\sigma f \Big) a_0 |\frac{d\sigma}{\sigma} dy'' d\mu''|^\half = 0, \ f \in
\CI(\Lpl),
\ilabel{left-eqn-0}
\end{equation}
near $T_+(\lam)$, or 
\begin{equation} 
\Big( -i\big({\mathcal L}_{V_l} 
- n \tau' \big) + 
\sigma f \Big) a_0 |\frac{d\sigma}{\sigma} \frac{ds}{s} dy'' d\muh''|^\half
= 0, \ f \in \CI(\Lpl),
\ilabel{left-eqn}
\end{equation}
away from $T_+(\lam)$,
which gives an equation for $a_0$ of the form
\begin{equation}
-i\tau' \big( \pa_\sigma  + f \big) (\sigma^{-n/2}a_0) = 0,
\quad f \in \CI(\Lpl).
\ilabel{eq-Vl}\end{equation}
This shows that $\sigma^{-n/2} a_0$ is smooth across $\sigma = 0$. 

To show regularity near $\Lsl$, we use the fact that the symbol satisfies
both the right and left transport equation. We take the sum of the transport
equations that obtain when we use the total boundary defining function $x'$
for $H_l$, and $x''$ for $H_r$. The right transport operator with respect to
$x''$ takes the form
$$
-i\big({\mathcal L}_{V_r} - n \tau'' \big) + p_{\sub} 
$$
However, by \eqref{change} the symbol written in terms of $x''$ is equal to
$(x''/x')^{-1/2-n/2}$ times the symbol written in terms
of $x'$. Since we are writing the symbol in terms of $x'$, we must include a
factor $\sigma^{-1/2-n/2}$ to be consistent with \eqref{left-eqn}. 
This gives an equation for $a_0$ of the form
\begin{equation} 
\Big( -i\big({\mathcal L}_{V_r} 
- n \tau'' \big) + 
p_{\sub}(y'',\mu'',\tau'') \Big)  
\Big( \sigma^{-1/2-n/2} a_0 |\frac{d\sigma}{\sigma} \frac{ds}{s} dy''
d\muh''|^\half \Big) = 0.
\ilabel{right-eqn}
\end{equation}
In view of the term $-2\tau'' \sigma \pa_\sigma$ in
the formula \eqref{Vr} for $V_r$, and since $p_{\sub}$ vanishes at $s=0$ since
$\mu'' = 0$ there, we get an equation for $a_0$ 
\begin{equation}
\big( V_r + \tau'' + s f' \big) a_0 = 0 \quad f' \in \CI(\Lplh).
\ilabel{eq-Vr}
\end{equation}
Combining with the left transport equation gives an equation which, using
\eqref{VlVr} and the fact that $\tau' = \tau'' = -\lam$ at $s=0$ takes the form
$$
2\lam \big( s \pa_s - \frac{n-1}{2} + sW + s\tilde f \big) a_0 = 0,
$$
where $W$ is tangent to the boundary of $\Lplh$ and $\tilde f$ is smooth on
$\Lplh$. This may be written
$$
\big( \pa_s + W + \tilde f \big) \big( s^{-(n-1)/2} a_0 \big) = 0,
$$
This together with \eqref{eq-Vl} shows that $a_0 \in \sigma^{n/2} s^{(n-1)/2}
\CI(\Lplh)$. 

\

It follows that there is a Legendre distribution $I^{-\frac1{2},
  \frac{n-2}{2}; \frac{n-1}{2}, \frac{n-1}{2}}(\Lpl,
\Lplh)$ which has the correct symbol of order $-1/2$ at $\Lpl$. Thus it
solves the equation
\begin{equation}
(H_l - \lam^2) Q - E(\lam) \in I^{\frac{3}{2}, \frac{n}{2};
 \frac{n+3}{2}, \frac{n-1}{2}}(\XXb, \Lpl, \Lsl; \sfOh),
\end{equation}
where $(n+3)/2$ is the order of
vanishing at lb and $(n-1)/2$ is the order of vanishing at rb. The order
of improvement at lb is two since not only is the Legendrian $G_1$ of
Definition~\ref{M-Leg-submfld} at lb characteristic for $H_l - \lam^2$, but
the transport 
operator for symbols of order $(n-1)/2$ vanishes, so we automatically get two
orders of improvement there. At rb however we can expect no improvement. 
As shown above, $Q$ will automatically satisfy the equation
$$
(H_r - \lam^2) Q - E_2(\lam) \in I^{\frac{3}{2}, \frac{n}{2};
 \frac{n-1}{2}, \frac{n+3}{2}}(\XXb, \Lpl, \Lsl; \sfOh).
$$
Let us assume by induction that we have a kernel which solves the left equation
above up to an error in 
\begin{equation}
I^{k+\frac{1}{2}, \frac{n}{2};
 \frac{n+3}{2}, \frac{n-1}{2}}(\XXb, \Lpl, \Lsl; \sfOh)
\ilabel{error-l-k}
\end{equation}
and hence the right equation up to an
error in 
\begin{equation}
I^{k+\frac{1}{2}, \frac{n}{2};
 \frac{n-1}{2}, \frac{n+3}{2}}(\XXb, \Lpl, \Lsl; \sfOh)
\end{equation}
We wish to improve this by one order
at $\Lpl$. To do this, we choose $Q_k \in I^{k-\frac1{2}, \frac{n-2}{2}; 
\frac{n-1}{2}, \frac{n-1}{2}}(\Lpl, \Lplh)$ to have the symbol of order
$k-1/2$ on $\Lpl$ 
which solves the left transport equation (and
therefore the right transport equation) on $\Lpl$. We need to investigate 
the regularity of this symbol to see if it extends to a Legendrian conic
pair. The argument is very analogous to the one above, but now we have error
terms on the right hand side. In the first region, after removing the
half-density factor, we get an equation of the form
\begin{equation}
-2i\tau' \big( \sigma \pa_\sigma + (-\frac{n}{2} + k) + \sigma f \big) a_k =
b_k, 
\quad f \in \CI(\Lpl).
\ilabel{eq-Vl-k}\end{equation}
The term $b_k$ comes from the error to be solved away. Since the error term
is of order $k+1/2$ at $\Lpl$ and order $(n+3)/2$ at lb, 
$b_k \in \sigma^{n/2-k+1}
\CI(\Lpl)$. This shows that $a_k \in \sigma^{n/2-k}
\CI(\Lpl)$, as desired. Similarly, in the second region, near the corner
$\lb \cap \Lsl$, by combining the vector fields $V_l + V_r$ we get an
equation of the form 
\begin{equation}
-2i\tau'' \Big( - s \pa_s + s W +( - \frac{n-1}{2} + k) + s \sigma f
\Big) a_k = b_k, \quad f \in \CI(\Lpl)
\ilabel{eq-Vr-k}\end{equation}
with $b_k$ again the error to be solved away. To calculate its order of
vanishing at $s=0$, consider the transport
equation for symbols of order $(n-2)/2$ at $\Lsl$. Noting that the
subprincipal symbols vanish identically on $\Lsl$, the left transport
operator is
$$ 
-i \big( \mathcal{L}_{V_l} -  \tau' \big) 
$$
whilst the right transport operator with respect to $x''$ is 
$$ 
-i \big( \mathcal{L}_{V_r} -  \tau'' \big) 
$$
To write this with respect to $x'$ we must conjugate by $\sigma$ (by
\eqref{change}). In view of the term $-2 \tau'' \sigma \pa_\sigma$, this
changes the operator to
$$ 
-i \big( \mathcal{L}_{V_r} +  \tau'' \big) .
$$
The sum of these two operators vanishes on $\Lsl$ so actually the right hand
side in  
\eqref{eq-Vr-k} comes from a term in $I^{k+\frac{1}{2}, \frac{n+2}{2};
 \frac{n+1}{2}, \frac{n-1}{2}}(\XXb, \Lpl, \Lsl; \sfOh)$. From
Proposition~\ref{ex-conic-2} we
see that $b_k \in s^{(n-1)/2-k+1} \sigma^{n/2-k} \CI(\Lplh)$, 
one power in $s$ better than might be expected. 
This shows that $a_k \in s^{(n-1)/2-k} \sigma^{n/2-k}
\CI(\Lplh)$, as desired. Therefore, one can find a Legendre conic pair with
symbol of order $k-1/2$ on $\Lpl$ equal to $a_k$ which solves away the error
term of order $k+1/2$ at $\Lpl$. 
This completes the inductive step. By asymptotically
summing these correction terms, we end up with an approximation $G_3(\lam)$ to the
resolvent kernel with an error $E_3(\lam)$ in $I^{\frac{n-2}{2};
\frac{n+3}{2}, \frac{n-1}{2}}(\XXb, \Lsl; \sfOh)$. That is, we have solved
away the scattering wavefront set 
of the error term at $\Lpl$ completely. 

\subsection{Solving away outgoing error}\label{outgoing-error}
In the last step of the construction of the parametrix, we solve away the
error to infinite order at bf and lb. We begin by considering the expansion
at rb. By construction, the parametrix $G_3(\lam)$ constructed so far has an
expansion at rb
$$
G_3(\lam) \sim e^{i\lam/x''} {x''}^{(n-1)/2} \sum_{j \geq 0} {x''}^j
g_j(z',y'') \cdot \nu' \cdot \nu'',
$$
where $g_j(z',y'') \cdot \nu' \in I^{-n/4-j, n/4 - 1/2-j}(G_{y''}(\lam),
\Gs(\lam))(X; \scOh)$,
$G_{y''}(\lam)$ is the fibre Legendrian of Definition~\ref{M-Leg-submfld} and
$y''$ is regarded as a smooth parameter. The factors $\nu'$, $\nu''$ are the
Riemannian half-density factors on $X$ lifted to $\XXb$ via the left and
right projections, respectively. We will ignore the half-density factors
from here on; since $\Lap(a \cdot \nu' \cdot \nu') = \Lap(a) \cdot \nu'
\cdot \nu'$, this only has the effect of changing $H = \Lap + P$ to $\Lap +
P'$ for some $P'$ with the same properties as $P$. 

The error term after applying $H - \lam^2$ to $G_3(\lam)$ has the form
\begin{equation}
E_3(\lam) \sim e^{i\lam/x''} {x''}^{(n-1)/2} \sum_{j \geq 0} {x''}^j
e_j(z',y''),
\ilabel{error3}\end{equation}
where $e_j \in {x'}^{(n+1)/2-j} e^{i\lam/x'} \CI(X \times Y)$. Again we regard
$y''$ as a parameter. Thus we have
\begin{equation}
(H - \lam^2) g_j = e_j  
\ilabel{gj}\end{equation}

Consider the problem of solving away errors of the
form $e_j$, to infinite order at bf (of course we cannot solve the errors
away exactly without begin able 
to solve $(H - \lam^2) u = f$ exactly, which we cannot do until we have
constructed the resolvent kernel!). 
If we apply $H - \lam^2$ to a series of the form 
\begin{equation}
e^{i\lam/x} x^{(n-1)/2 - k}\sum_{j \geq 0} x^j b_j, \quad b_j \in \CI(X),
\ilabel{series1}\end{equation}
we get a series of the form
\begin{equation} 
e^{i\lam/x} x^{(n+1)/2 - k} \sum_{j \geq 0} x^j c_j, \quad c_j \in \CI(X),
\quad c_0 = 2i\lam k b_0.
\ilabel{series2}\end{equation}
Thus, we can add to $g_j$ a series of the form \eqref{series1} to solve away
the powers greater than $(n+1)/2$, but the power $(n+1)/2$ presents a
problem (without introducing logarithmic terms), because of the vanishing of
$k$ in \eqref{series2}, unless the 
coefficient of the $x^{(n+1)/2}$ term happens to be zero. We need the
following results.

\begin{lemma}\label{lemma:one-order-improvement}
If $g \in  e^{i\lam/x} x^{(n-1)/2 - k} \CI(X)$, $k = 1, 2, \dots$, satisfies 
$$
(H - \lam^2) g = x^{(n+1)/2} e^{i\lam/x} \CI,
$$
then 
$$
(H - \lam^2) g = x^{(n+3)/2} e^{i\lam/x} \CI.
$$
\end{lemma}

\begin{proof} It follows inductively using \eqref{series2} that the
coefficient of order $(n+1)/2 - l$ vanishes, for $l = k, k-1,
\dots$. Thus, actually $g \in e^{i\lam/x} x^{(n-1)/2} \CI(X)$. Then
\eqref{series2} 
shows that the next coefficient also vanishes. 
\end{proof}

\begin{cor} The same result holds if the condition $g \in  e^{i\lam/x}
x^{(n-1)/2 - k} 
\CI(X)$ is replaced by $g \in I^{p, (n-1)/2 - n/4 - k}(K, \Gs)$ for any
$p$ and any Legendre conic pair $(K, \Gs)$.
\end{cor}

\begin{proof} Apply the above argument to the symbol at $\Gs$.
\end{proof}

Thus, for each $j$, we can modify $g_j$ by a series of the form
\eqref{series1} until the error term is of the form ${x'}^{(n+1)/2}
e^{i\lam/x'} \CI(X \times Y)$. Then applying the Corollary to $g_j$,
we find that unsolvable term of order ${x'}^{(n+1)/2}$ vanishes. 
Therefore, we can solve away the $e_j$ to infinite order
at bf. Thus, we may assume that our error in $E_3(\lam)$ vanishes to infinite
order at the corner bf $\cap \rb$.

Next, we solve the error away at $\Ls$. This involves solving the transport
equation 
\begin{equation}
i\lam \Big( \sigma \pa_{\sigma} + (\frac{1}{2} + j)
\Big) a_j  = b_j.
\end{equation}
The equation for $a_0$ then is
$$
\Big( \sigma \pa_{\sigma} + \frac{1}{2} \Big) a_0 = b_0,
$$
and $b_0$ is rapidly decreasing at rb and is in $\sigma^{3/2} \CI(\Lsl)$ at
lb. There is a unique solution which is rapidly decreasing at rb and in
$\sigma^{-1/2} \CI$ at lb. We can thus find a correction term which reduces
the error to $I^{(n+2)/2; (n+3)/2, (n-1)/2}(\Lsl)$, 
with infinite order vanishing at bf $\cap \rb$. 
Inductively, assume that we have reduced the error to $I^{n/2 + k; (n+3)/2,
(n-1)/2}(\Lsl)$, with infinite order vanishing at bf $\cap \rb$. 
The transport equation for $a_k$ is then
$$
\Big( \sigma \pa_{\sigma} + \frac{1}{2} + k \Big) a_k = b_k,
$$
where inductively, $b_k$ is rapidly decreasing at rb and is in $\sigma^{1/2 - k}
\CI(\Lsl)$ at lb. There is a unique solution rapidly decreasing at
rb and in $\sigma^{-1/2 + k} \CI(\Lsl)$ at lb. A Legendre distribution in 
$I^{n/2 + k - 1, (n-1)/2, (n-1)/2}(\Lsl)$ with $a_j$ as symbol then reduces
the error to $I^{n/2 + k + 1; (n+3)/2,
(n-1)/2}(\Lsl)$, with infinite order vanishing at bf $\cap \rb$, so this
completes the inductive step. Taking an asymptotic sum of such correction
terms yields a parametrix $G_4(\lam) = G(\lam)$ leaving an error which is the
sum of a term supported away from rb 
of the form 
$$
e^{i\lam/x'} {x'}^{(n+3)/2} a(y', z'')
$$
with $a$ smooth and rapidly decreasing at bf, plus a term supported away
from lb of the form
$$
e^{i\lam/x''} {x''}^{(n-1)/2} b(y'', z')
$$
at rb, with $b$ rapidly decreasing at bf. The error at lb can be solved away
using \eqref{series1}-\eqref{series2}, leaving an error term $E_4(\lam)$
which can be expressed on the
blown-down space $X^2$ as 
$$
E_+(\lam) = e^{i\lam/x''} {x''}^{(n-1)/2} b(z', z'')
$$
with $b$ smooth on $X^2$ and rapidly decreasing at $x' = 0$. Such an error
term is compact on the weighted $L^2$ space $x^l L^2(X)$ for any $l > 1/2$
(where $L^2$ is taken with respect to the metric density). 
This completes the construction of the parametrix.

\section{Long range case}
The case of long range metrics and long range perturbations, 
$P\in x\Diffsc^{1}(X)$, requires only
minor modifications in the parametrix construction until the last step
(removing the outgoing error). In particular,
there is no change in the construction of the pseudodifferential
approximation.
In the intersecting Legendrian construction, as well as solving the
transport equations on $L^\circ(\lambda)$, the only difference is in the
structure of the subprincipal symbol, which no longer obeys
Lemma~\ref{lem:sub-pr}. Thus, the arguments in section~\ref{smoothness} and section~\ref{outgoing-error} have
to be modified. Let $q$ denote the boundary subprincipal symbol of
$H$. Notice that in the gravitational long range case, $q$ is a constant,
but in the general long range case, $q$ is an arbitrary smooth function of $y$
which is a quadratic on each fibre of $K$ over $Y$. 
Let $q_l$ and $q_r$ denote the lift of $q$ to $\XXb$ via the left,
respectively right, projection. 

Let us now discuss the necessary modifications to sections~\ref{smoothness}
and \ref{outgoing-error}. Equation \eqref{left-eqn-0} becomes
\begin{equation} 
-i\big({\mathcal L}_{V_l} 
- n \tau' +iq_l\big)  a_0 |\frac{d\sigma}{\sigma} dy'' d\mu''|^\half = 0,
\end{equation}
near $T_+(\lam)$. Thus,
\eqref{eq-Vl} is replaced by
\begin{equation}
-i\tau' \big( \pa_\sigma  + f \big) (\sigma^{-n/2
-i\frac{q_l}{2\lambda}}a_0) = 0,
\quad f \in \CI(\Lpl).
\end{equation}
Thus, now we conclude that in this region
$a_0$ is of the form
$$
\sigma^{\frac{n}{2}+i\frac{q_l}{2\lambda}}\CI(\Lpl).
$$

Next, in the second region, at the corner $\lf\cap L^\sharp(\lam)$,
the right transport equation \eqref{eq-Vr} becomes
\begin{equation}
\big( V_r + \tau'' +iq_r+ s f' \big) a_0 = 0
\quad f' \in \CI(\Lplh).
\end{equation}
Adding this to the left transport equation yields
$$
2\lam \big( s \pa_s - \frac{n-1}{2}-i\frac{q_l}{2\lam}-i\frac{q_r}{2\lam}
 + sW + s\tilde f \big) a_0 = 0,
$$
which now gives that $a_0$ is of the form
$$
\sigma^{\frac{n}{2}+i\frac{q_l}{2\lambda}}
s^{\frac{n-1}{2}+i\frac{q_l + q_r}
{2\lambda}}\CI(\Lplh).
$$

Combining this with the other similar results at $\rf$ and the interior
of $L^\sharp(\lam)$, we deduce that
$$
a_0\in\rho_{\lb}^{n/2+i\frac{q_l}{2\lambda}}
\rho_{\rb}^{n/2+i\frac{q_r}{2\lambda}}
\rho_{\#}^{(n-1)/2+i\frac{q_l + q_r}{2\lambda}}
\CI(\Lplh;
\Omega_b^{1/2} \otimes S^{[-1/2]}).
$$ 

In the general long range case, 
the dependence of $q$ on $y$, and its appearance 
in the exponent of the boundary defining functions
$\rho_{\lb}$, etc.,
means that differential operators from the left
factor, acting on a Legendre function with principal symbol $a_0$,
introduce logarithmic terms. For example, in a neighborhood
of $\lb$ in $\Lpl$ the error term $b_k$
in \eqref{eq-Vl-k} for $k=1$ will
take the form
$$
b_1=\sigma^{\frac{n}{2}+i\frac{q_l}{2\lambda}}
((\log\sigma)^2 c_2+\log\sigma\, c_1+c_0),\quad c_j\in\Cinf(\Lpl),\ j=0,1,2.
$$
Then the transport equation for $a_1$ takes the form
\begin{equation} 
\Big( -i\big({\mathcal L}_{V_l} 
+(- n+2) \tau' +iq_l\big) +
\sigma f \Big) a_1 |\frac{d\sigma}{\sigma} dy'' d\mu''|^\half = 0, \ f \in
\CI(\Lpl),
\end{equation}
\begin{equation}
-i\tau' \big( \sigma \pa_\sigma  -\frac{n}{2}+1 -iq_l/(2\lambda)
+ \sigma f \big) a_1 = b_1,
\quad f \in \CI(\Lpl).
\end{equation}
Hence, near $\lb$ but away from $L^\sharp(\lam)$,
$a_1$ will take the form
$$
a_1=\sigma^{\frac{n}{2}-1+i\frac{q_l}{2\lambda}}
((\log\sigma)^2 c'_2+\log\sigma\, c'_1+c'_0),
\quad c'_j\in\Cinf(\Lpl).
$$

A similar discussion works at the other boundary faces of $\Lh^+(\lam)$,
with up to quadratic factors in each of $\log \rho_{\lb}$,
$\log \rho_{\rb}$, $\log_{\rho_\sharp}$, and can be repeated (with
progressively higher powers of logarithms) for all $a_k$'s.

Since the most important long-range case is the gravitational case
where the subprincipal symbol is constant, 
and since it makes the discussion more transparent,
in what follows we make the assumption that
$$
\Lap \text{ and } P \ \text{are of long range gravitational type},
$$
which implies that $q$ is constant. Let
$$
\alpha = \frac{q}{2\lam}.
$$
The point is that in this case the powers of $\rho_{\lb}$, etc., above
are constant, thus no logarithmic factors arise when we apply $H-\lam^2$
to such Legendre functions.
Then
$$
a_k\in\rho_{\lb}^{n/2+i\alpha}
\rho_{\rb}^{n/2+i\alpha-k}
\rho_{\#}^{(n-1)/2+2i\alpha-k}
\CI(\Lplh;
\Omega_b^{1/2} \otimes S^{[-1/2]}),
$$
and asymptotic summation gives an outgoing error
$$
E_+(\lam) \in e^{i\lam/x'}e^{i\lam/x''}{x''}^{(n-1)/2+i\alpha}
{x'}^{(n+1)/2+i\alpha} \CI(\XXb;
\scOh).
$$
Since $\alpha$ is a constant, \eqref{series1}-\eqref{series2} are
still true ($k$ need
not be an integer; it suffices that it is a constant), except that
now
$$
c_0=(2i\lam k+2\lam\alpha)b_0 .
$$
Since now we are taking $k=l+i\alpha$,
where $l$ is an integer, we can solve away the series, provided that
the coefficient of the $l=0$ term vanishes, which is assured just as in
Lemma~\ref{lemma:one-order-improvement}. The rest of the argument
requires only similar modifications as compared to the short-range case,
so we conclude, as there, that we can modify the parametrix to obtain
an error term of the form
$$
E_4(\lam)=e^{i\lambda/x''}(x'')^{(n-1)/2+i\alpha}b(z',z'')
$$
with $b$ smooth on $X^2$ and rapidly decreasing at $x'=0$.

\section{Resolvent from parametrix}

In the previous two sections, we constructed a parametrix $G(\lam)$ for $\Rt(\lam)$ 
which satisfies
$$
(H_l - \lam^2) G(\lam) = K_{\Id} + E(\lam),
$$
where $E(\lam)$ has a kernel which is of the form $e^{i\lam/x''} {x''}^{(n-1)/2} 
{x'}^\infty \CI(\XXb; \scOh)$. Thus, it is a Hilbert-Schmidt kernel on
$x^l L^2(X)$ for every $l > 1/2$, and in particular is compact. In fact,
we also see directly from its form that $E(\lam):x^l L^2(X)\to\dCinf(X)$
for $l>1/2$. Crude
estimates (such as Schur's Lemma)
show that $G(\lam)$ acts as a bounded operator from $x^l L^2$ to $x^{-l} L^2$
for large enough $l>1/2$;
more refined estimates, which we do not need here,
show that in fact this is true for any $l >
1/2$. Thus, the equation above becomes an operator equation
$$
(H - \lam^2) G(\lam) = \Id + E(\lam)
$$
from $x^l L^2$ to $x^{-l} L^2$. 

\subsection{Finite rank perturbation}
To correct $G(\lam)$ to the actual $\Rt(\lam)$,
we must solve away the error term $E(\lam)$. Thus, we would like $\Id + E(\lam)$
to be invertible. However, this is certainly not
necessarily the case as things stand; if for example we modified $G(\lam)$ 
by subtracting from it the rank one operator $G(\lam)(\phi) \langle \phi,
\cdot \rangle$, for some $\phi \in \CIdot(X; \scOh)$, then the modified $G(\lam)$
would be microlocally indistinguishable from the old one, but would
annihilate $\phi$, so the modified $\Id + E(\lam)$ would not be invertible.  

Since $\Id + E(\lam)$ is compact, it has a null space and cokernel of the same
finite dimension $N$. 
To make $\Id + E(\lam)$ invertible, we try to correct $G(\lam)$ by adding to it a
finite rank term 
\begin{equation}
\sum_{i=1}^N \phi_i \langle x^{2l} \psi_i, \cdot \rangle .
\ilabel{finiterank}\end{equation}
Here $\langle, \ \rangle$ denotes $L^2$ pairing, $\psi_i$ should lie in
$L^2$,
and the factor of $x^{2l}$ is included to ensure that it acts
on $x^l L^2$. We require that $\phi_i$ are in $x^{-l} L^2$ so that
\eqref{finiterank} is bounded from $x^{l} L^2 \to x^{-l} L^2$. We wish to
choose $\phi_i$ and $\psi_i$ so that 
$$
\Id + E(\lam) + \sum_{i=1}^N ((H - \lam^2) \phi_i) \langle x^{2l} \psi_i, \cdot
\rangle 
$$
is invertible. This is possible if we can choose $x^{l} \psi_i$ to span the
null space of $\Id + E(\lam)$ and $(H - \lam ^2)\phi_i$ to span a subspace
supplementary to 
the range. Note that if $(\Id+E(\lam))u=0$ and $u\in x^l L^2$, then
$u=-E(\lam)u$, so the mapping properties of $E(\lam)$ imply that
$u\in\dCinf(X)$. Thus, we automatically have
$\psi_i\in\dCinf(X)$ above. To proceed, we need the following lemma. 

\begin{lemma}\label{density}
Let $l>1/2$. Then the image of $H - \lam^2$ on
the sum of $\dCinf(X)$ and the range of $G(\lam)$ applied to $\dCinf(X)$
is dense in $x^l L^2$.  
\end{lemma}

\begin{rem}
Note that $(H-\lam^2)G(\lam)g=(\Id+E(\lam))g\in\dCinf(X)$ if $g\in\dCinf(X)$,
and for $u\in\dCinf(X)$, $(H-\lam^2)u\in\dCinf(X)$, so the image of
$H-\lam^2$ on the space in the statement of the lemma is a subspace
of $\dCinf(X)$.
\end{rem}

\begin{proof} To proceed, we give the proof for short range $H$; the proof
for long range $H$ requires only minor modifications. 

Let $\mathcal{M}$ be the subspace of $x^l L^2$ given by the 
image of $H - \lam^2$ on
the sum of $\dCinf(X)$ and the range of $G(\lam)$ applied to $\dCinf(X)$. 
If $\mathcal{M}$ is not dense, then there is a function $f \in x^l L^2$
orthogonal to 
$\mathcal{M}$. Since for
$u\in\dCinf(X)$ implies
$(H-\lam^2)u\in\mathcal{M}$, $f$ satisfies
\begin{equation}\begin{gathered}
\lang x^{-l} f, x^{-l} (H - \lam^2) u \rang = 0 \quad \forall \ u \in
\CIdot(X) \\
\Rightarrow \lang (H - \lam^2 ) x^{-2l} f, u \rang = 0 \quad \forall \ u \in
\CIdot(X) \\
\Rightarrow (H - \lam^2) h = 0, \quad h = x^{-2l} f. 
\end{gathered}\end{equation}
where we used that $H$ is symmetric on $\dCinf(X)$.
On the other hand, $G(\lam)$ maps $\dCinf(X) \to x^{(n-1)/2}
e^{i\lam/x} \CI(X)$, and for any
$g\in\dCinf(X)$, $(H-\lam^2)G(\lam)g=(\Id+E(\lam))g\in\dCinf(X)$, hence
$(\Id+E(\lam))g\in\mathcal{M}$. In addition, $E(\lam)^*$, with kernel
$E(\lam)^*(z',z'')
=\overline{E(\lam,z'',z')}$, maps $\dist(X)\to x^{(n-1)/2}e^{-i\lambda/x}
\CI(X)$,
so we have
\begin{equation}\begin{gathered}
\lang x^{-l} f, x^{-l} (H - \lam^2) G(\lam) g \rang =
\lang x^{-l} f, x^{-l} (\Id + E(\lam)) g \rang = 0 \quad \forall \ g \in
\CIdot(X) \\
\Rightarrow \lang (\Id + E(\lam)^* ) x^{-2l} f, g \rang = 0 \quad \forall \ g \in
\CIdot(X) \\
\Rightarrow (\Id + E(\lam)^*) h = 0. 
\end{gathered}\end{equation}
If $h = - E(\lam)^* h$, then $h$ has the form $x^{(n-1)/2} e^{-i\lam/x}
\Cinf(X)$, ie it is incoming. A standard argument then implies that $h
\equiv 0$. Let $h = x^{(n-1)/2} e^{-i\lam/x} h_0(y) + \tilde h$, where $\tilde
h \in x^{(n+1)/2} e^{-i\lam/x} \CI(X)$. Green's formula yields
\begin{equation}\begin{gathered}
0 = \int_X \overline h (H - \lam^2) h  - h (H - \lam^2) \overline h \\
= 2i\lam \lim_{\ep \to 0} \int\limits_{\{ x = \ep \}} \overline h x^2 \pa_x h  -
h x^2 
\pa_x \overline h 
= 2i\lam \int\limits_Y |h_0(y)|^2,
\end{gathered}\end{equation}
so $h_0 \equiv 0$. It then follows iteratively from \eqref{series1} and 
\eqref{series2} that the expansion of $h$
at the boundary of $X$ vanishes identically, that is, that $h \in
\CIdot(X)$. Finally a unique continuation theorem, see
e.g.\ \cite[Chapter~XVII]{Hor}, shows $h=0$ identically. 
This means
that $\mathcal{M}$ is indeed dense in $x^{l} L^2$.
\end{proof}

Thus, we can choose the $\phi_i \in x^{(n-1)/2}
e^{i\lam/x} \CI(X)$ above so that $(H - \lam^2) \phi_i\in\dCinf(X)$
span a supplementary subspace of range $\Id + E(\lam)$.
The modified parametrix then satisfies 
$$
(H  -\lam^2) G_5(\lam) = \Id + E_5(\lam),
$$
where $E_5(\lam)$ has the same form as $E(\lam)$ but in addition is 
invertible on $x^l L^2$ for all $l > 1/2$. 

\subsection{Resolvent}
 The inverse $\Id + S(\lam)$ of $\Id + E_5(\lam)$ is Hilbert-Schmidt on $x^l
 L^2$ since this is true of $E_5(\lam)$. Moreover, since 
$$
S(\lam) = -E_5(\lam) + E_5(\lam)^2 + E_5(\lam) S(\lam) E_5(\lam),
$$
it is easy to see that also $S(\lam) \in e^{i\lam/x''} {x''}^{(n-1)/2} 
{x'}^\infty \CI(\XXb; \scOh)$. Our solution for the kernel $\Rt(\lam)$ is
then
$$
\Rt(\lam) = G_5(\lam) ( \Id + S(\lam)).
$$
It is not hard to show that $G_5(\lam) S(\lam)$ has the form
$${x'}^{(n-1)/2}
{x''}^{(n-1)/2} e^{i\lam/x'} e^{i\lam/x''} \CI(\XXb; \scOh),$$
so $\Rt(\lam)$
has the desired microlocal regularity \eqref{seek}. 

{\it Remark. } Lemma~\ref{density} directly shows the absence of positive
eigenvalues of $H$. Suppose that $(H - \lam^2) u = 0$ and $u \in x^s H^2(X)$
for some $s > -1/2$. This would certainly be the case of an eigenfunction
since $H$ has elliptic interior symbol, so $u$ would lie in $H^k(X)$ for
every $k$. We need to show that
for all functions $g\in x^{(n-1)/2}e^{i\lam/x}
\Cinf(X)$ the equation
\begin{equation}
\lang (H - \lam^2) u, g \rang = \lang u, (H - \lam^2) g \rang = 0
\ilabel{ef-symmetry}
\end{equation}
holds. Indeed, this implies
that $u$ is $L^2$-orthogonal to the image of $H - \lam^2$
acting on
$x^{(n-1)/2}e^{i\lam/x}
\Cinf(X)$,
or equivalently that $x^{2l} u \in x^l L^2$ is orthogonal in
$x^l L^2$ to the image of $H - \lam^2$ acting on $x^{(n-1)/2}e^{i\lam/x}
\Cinf(X)$ in $x^l L^2$. Then
Lemma~\ref{density} shows that $u \equiv 0$.

To deduce \eqref{ef-symmetry} for $g\in x^{(n-1)/2}e^{i\lam/x}
\Cinf(X)$,
let $\phi\in\Cinf(\RR)$, identically $1$ on $[1,\infty)$, $0$ near the
origin.
Then
\begin{equation}
\begin{split}
\lang (H-\lam^2)u, g\rang&=\lim_{t\to 0}\lang (H-\lam^2)u, \phi(x/t)g\rang
=\lim_{t\to 0} \lang u,(H-\lam^2)\phi(x/t)g\rang\\
&=\lim_{t\to 0} \lang u,\phi(x/t)(H-\lam^2)g\rang
+\lim_{t\to 0} \lang u, [H,\phi(x/t)]g\rang\\
&=\lang u,(H-\lam^2)g\rang+\lim_{t\to 0} \lang u, [H,\phi(x/t)]g\rang.
\end{split}\end{equation}
Note that $[H,\phi(x/t)]$ is uniformly bounded (i.e.\ with bounds
independent of $t$)
as a map $x^l H^1\to x^{l+1}L^2$, and in fact $[H,\phi(x/t)]\to 0$ strongly
as $t\to 0$. Applying this with $l=-1/2-\ep$, $\ep>0$ sufficiently small,
we see that the last term goes to $0$ as $t\to 0$, proving
\eqref{ef-symmetry}.

\subsection{Analytic continuation}
It is not hard to show that the kernel $\Rt(\lam)$ constructed
above continues analytically (as a distribution on $\XXsc$)
into $\Imag \lam \geq 0, \Rea \lam > 0$. We
complete the proof of Theorem~\ref{main-result} by showing that this 
analytic continuation coincides with the outgoing resolvent, $R(\lam^2)$,
for $\Imag \lam > 0$. 

Our parametrix is defined as an asymptotic sum of
symbols, which is really a sum with cutoff functions inserted
(see \cite[Proposition~18.1.3]{Hor} for an explicit construction). The
cutoffs
depend on $C^k$ norms of a finite number of symbols and ensure that the
sum
converges in $C^k$ for all $k$. If the symbols are
holomorphic in $\lam$
then the $C^k$ norms may be taken uniform on compact subsets of
$\lam$. Since holomorphy is preserved under uniform limits, we need only
check that the phase and symbols analytically continue in some explicit
parametrization of the Legendrians. 

It is standard that the pseudodifferential approximation $G_1(\lam)$
analytically continues. Blowing down sf, we solve away
the error as an intersecting Legendrian,
see Section~\ref{sec:intersecting-Leg-constr}. 
Let $\phi$ be a local parametrization of the Legendrian $L(1)$ near $L(1)
\cap N^*\diag_\bl$. Then it is easily checked that the phase function
$$
k \phi + s(\lam - k)
$$
locally parametrizes $(N^*\diag_\bl, L^+(\lam))$. Since the variable $s$
takes
nonnegative values, the function $e^{i(k\phi + s(\lam - k))/x}$ 
continues to $\Imag \lam > 0$.

Away from sf, the value of $\tau$ is strictly negative on the
Legendrian, and so the phase is of the form
$$
e^{i\lam \phi/x}
$$
where $\phi$ is positive on the Legendrian, independent of $\lam$.
By restricting the support of
the symbol sufficiently, therefore, we may assume that $\phi$ is positive
everywhere on the support of the integral. Thus this also analytically 
continues to the upper half plane with uniform bounds.

The symbols are defined by iteratively solving transport equations
of the form
$$
\Big( -i \mathcal{L}_{H_p} -i \big(\half + m - \frac{N}{4} \big) \frac{
\pa
  p}{\pa \tau} + p_{\sub} \Big) a_j=b_j,
$$
where $b_0=0$.
These equations
are solved along $L^\circ(\lam)$, i.e.\ if we consider amplitudes
in an explicit parameterization of the Legendrian, then along
the critical submanifold $C_\phi=\{(0,y,u):\ d_u\phi=0\}$, where
$\lambda\phi/x$ is the phase function as above. Note that $C_\phi$
is independent of $\lam$, and it is identified with $L^\circ(\lam)$
via the map $C_\phi\ni(0,y,u)\mapsto(0,y,d_{(x,y)}(\lam\phi/x)|_{(0,y,u)})$.
Along $C_\phi$ the transport equation
becomes an ODE whose coefficients depend on $\lam$ polynomially, since the
only $\lam$ dependence of the coefficients arises from this identification
map,
and $H_p$, $\frac{ \pa  p}{\pa \tau}$, $p_{\sub}$
are polynomial in the fiber variables.
Thus, the solution $a_j$ of the transport equation,
as a function on $D\times C_\phi$, $D$ a neighborhood of the positive
real axis in $\CC=\CC_\lam$, is
holomorphic (in $\lam$), provided $b_j$ is (here we
identify $C_\phi$ with $L^\circ(\lam)$). Note that
the $b_j$'s arise because solving the transport equations only
guarantees that the `error term' $E_3(\lam)$, arising from the application
of $H-\lam^2$ to $R_3(\lam)$, is one order lower than expected, so for each
$\lam$,
$b_j|_{C_\phi}$
depends on $a_i$, $i<j$ near $C_\phi$, and not just on $a_i|_{C_\phi}$.
(In fact, $b_j|_{C_\phi}$ depends on a finite number of terms of the
Taylor series of $a_i$, $i<j$, at $C_\phi$.)
To ensure that the
$b_j$ are holomorphic in $\lam$, we define the $a_i$ {\em near} $D\times
C_\phi$,
rather than {\em at} $D\times C_\phi$, e.g. by introducing a local product
decomposition $C_\phi\times U$, $U\subset\RR^k$, of the parameter space
near $C_\phi$, and pull-back the $a_i$, first defined on $D\times C_\phi$,
by the projection. Then, having constructed $a_i$, $i<j$, $b_j$ will
be holomorphic in $\lam$ near, hence at, $D\times C_\phi$, so $a_j$ is also
holomorphic at $D\times C_\phi$, hence it extends to be holomorphic near
$D\times C_\phi$. If we express the amplitudes $a_j$ with respect to
a different parameterization of $L^\circ(\lam)$, which
is still of the form $\lambda\phit/x$, then the new amplitudes $\at_j$
will
still be holomorphic functions of $\lam$, so holomorphy is preserved
at the overlap of parameterizations of different parts of $L^\circ(\lam)$.
This completes the inductive argument.

Therefore, our parametrix constructed above may be assumed holomorphic in
some set $B(\ep, \lam_0) \cap \{ \Imag \lam \geq 0 \}$, for some $\lam_0 > 0$.
It
is
easy to see that for non-real $\lam$, the parametrix is in the small
calculus, since the positivity of $\phi$ implies that the exponent
of $e^{i\lam\phi/x}$ has negative real part, and is therefore rapidly
decreasing at bf, lb and rb. 
The finite rank correction may be taken independent of $\lam$ if we
chose $\ep>0$ sufficiently small. Then, 
we have
$$
(\Lap  -\lam^2) G_5(\lam) = \Id + E_5(\lam), \quad \Imag \lam \geq 0, \ \Rea \lam > 0,
$$
where all terms are holomorphic in some small open set as above, 
$E_5(\lam)$ is invertible on $x^l
L^2$ for all $l > 1/2$, and off the real axis, $G_5(\lam)$ and $E_5(\lam)$ are in
the small calculus. Define $\Id + S(\lam)$ to be the inverse of $\Id +
E_5(\lam)$ on $x^l L^2$ for some fixed $l$. By the symbolic functional
calculus \cite{HV1}, for $\Imag \lam > 0$,
$S$ is a family of scattering pseudodifferential operators which is
clearly
holomorphic. Then $\Rt(\lam) \equiv G_5(\lam)(\Id + S(\lam))$
satisfies $(H - \lam^2)\Rt(\lam) = \Id$ on $x^l L^2$. But by
self-adjointness,
and the symbolic functional calculus, 
for $\Imag \lam > 0$, $(H - \lam^2)$ has a pseudodifferential inverse on
$L^2$. Since $R(\lam)$ is a bounded operator on $L^2$ for $\Imag \lam > 0$ it must
be the 
inverse. Therefore we have shown the inverse on the real axis constructed
above continues as a Schwartz kernel to the upper half plane and agrees
with the resolvent there. This completes the proof of
Theorem~\ref{main-result}. 

\begin{rem}
The only place where we used that $\Imag \lam\geq 0$ is to make our
parametrix act on, and its error compact on, weighted Sobolev
spaces. Namely, in the last step
of the construction, i.e.\ when we add a finite rank perturbation to
remove the error $E(\lam)$, we need $E(\lam)$ to be a compact
operator on $x^l L^2$ for $l>1/2$. However, the kernel of $E(\lam)$ is of
the form
$e^{i\lam/x''} {x''}^{(n-1)/2} 
{x'}^\infty \CI(\XXb; \scOh)$, and for $\Imag \lam<0$ the real part of the
exponent is positive, so the kernel of $E(\lam)$ is not even a tempered
distribution on $\XXb$. In particular, it does not even map $\dCinf(X)$
to $\dist(X)$. The same statement holds for $G(\lam)$ as well.
\end{rem}

\bibliographystyle{plain}
\bibliography{sm}
\end{document}